\documentclass[12pt,epic,eepic,psfig,epsfig]{article}
\usepackage{epic,eepic,epsfig,amssymb}
\def\endpf{{\ \hfill\hbox{\vrule width1.0ex height1.0ex}\parfillskip 0pt}}
\newenvironment{proof}{\noindent{\bf Proof:}}{\endpf}

\newenvironment{proof11}{\noindent{\bf Proof of Theorem 1S:}}{\endpf}

\newenvironment{proof1CP}{\noindent{\bf Proof of Theorem 1CP:  }}{\endpf}
\def\Pr{{\rm Pr}}
\def\x{\hbox{\b{\hbox{$x$}}}} 
\def\t{\hbox{\b{\hbox{$\theta$}}}} 

\begin{document}
\newcommand{\asinglelineskip}{\setlength{\baselineskip}{1.0\baselineskip
}}
\newcommand{\adoublelineskip}{\setlength{\baselineskip}{1.66\baselineskip}}
\newenvironment{exenv}{\asinglelineskip}{\adoublelineskip}

\title{{\bf A GEOMETRIC FORMULATION \\
OF FIDUCIAL PROBABILITY}
          \author{Paul Gunther
      \thanks{4515 38th Street NW, Washington, DC 20016, 
     e-mail: \{\rm igunther@rcn.com\}
     }}
                         \date{\vspace*{-50pt}}
\\
}
                         \maketitle

\abstract{The geometric formulation of fiducial 
probability employed in this paper is an  improvement over the usual
pivotal quantity formulation.
For a single parameter and single observation, the new formulation is based
on the geometric  properties of an ordinary two variable function and its
surface representation.

The following theorem is proved: 
A fiducial distribution for the continuous parameter $\theta$ exists if and
only 
if (i) the continuous random probability distributions of $x$ for different
$\theta$'s are non-intersecting, and (ii) the random
distributions are complete, i.e. at the extreme values of $\theta$ the
limiting probability distributions are zero and one for all $x$. 

The proof yields also  a complete characterization of random distributions
that lead to fiducial distributions.

The paper also treats intersecting distributions and
non-intersecting incomplete distributions. The latter, which are frequently
encountered in a null hypothesis, are shown to be associated
with intersecting ``composite'' distributions.

An appendix compares the pivotal and geometric formulations. 

\adoublelineskip

\section*{Introduction}

Fiducial  probability was introduced by Fisher in 1930 [2].
Since the 1960's, however, interest and research in the subject has
practically ceased and fiducial probability for the most part rejected. 
(An exception is the recent paper by Hannig [5] on ``generalized
fiducial inference''.) 

The approach adopted by Fisher, and generally accepted by researchers, 
can be characterized as a {\em pivotal quantity} (PQ) formulation.
The present paper presents an alternative geometric formulation for
fiducial probability that is based on the properties of any function of
two variables and its surface representation.

Appendix A, which compares the PQ and geometric formulations, 
also identifies flaws in the PQ formulation. To better understand the
problem it will be useful to briefly summarize the analysis as well as the
considerations that led to the PQ formulation.

Early on it was seen by Fisher that for translation-scale parameters a
pivoting quantity yielded both a random distribution and a fiducial
distribution.
Apparently the PQ was viewed as a necessary ingredient for obtaining a
fiducial distribution. 
Fisher also presented this result in the form of a general equation
relating
fiducial and random distributions, that was applicable to all
parameters.
The equation, which actually defines what we call the geometric
formulation, reflects the fact that the PQ formulation for
translation-scale parameters is also a  geometric formulation.
The PQ property was apparently so convincing, however, that there was no
need to pursue any alternative formulation.

A natural {\em pivoting extension} that provides a generalization to
all parameters was obtained by defining a general PQ by the equation
$F(x,\theta)=$constant. Pivoting consisted of inverting a monotone
random probability function to  obtain an equivalent confidence limit
function. 
Careful scrutiny of this procedure in  Appendix A shows, however, that this
generalization is neither  pivoting nor an
extension. 

The flaw in the PQ formulation may very well have contributed to the
subsequent decline of fiducial probability research. 

The benefits of the geometric formulation replacement 
are demonstrated by its achievements, such as the non-intersection condition
for existence of fiducial distributions. Also obtained is the solution for 
the fiducial distribution and associated confidence limits for multiple
observations (the derivation presented in Appendix A) that is applicable to
all parameters. This problem 
was  previously unsolved, ostensibly because it requires a
(geometric formulated) fiducial argument rather than the random probability
argument that sufficed for a single observation.  This result  also
demonstrates the necessity
for including fiducial theory as an integral part of
parametric probability theory. 

The new geometric formulation may well lead to renewed interest and
research into fiducial probability, especially in  areas of multiple
parameters and discrete parameters.

Except for Appendix A, the paper is devoted entirely to a non-intersection
existence 
theorem and its consequences. Section 1 illustrates the
determination of the fiducial distribution using the geometric formulation.
Section 2 gives an example of intersecting distributions that is useful
for the existence proof in Section 3. The proof includes a detailed
analysis of ``touching"  distributions, and leads also to a complete
characterization of non-intersecting  configurations.
Section 4 discusses the three dimensional fiducial surface and associated
geometry that underlies the geometric  formulation. 
Section 5 analyzes intersecting ``composite" distributions, and its
evolution into the exceptional case of non-intersecting incomplete fiducial
distributions, which arise in many null hypotheses.

\bigskip

\section{Geometric formulation of the fiducial
distribution}\label{section1}

This section uses a plane representation  of the two variable function
that forms the basis of the geometric formulation.
(Section 4 provides a three dimensional surface representation.)

\bigskip
\noindent{\bf Definition CFM}  
A {\em continuous fiducial model} for the single 
observation $x$ and
the single parameter $\theta$  satisfies 
the following three conditions:

\noindent (a) the domains of $x$ and $\theta$ are (finite or infinite)
intervals;

\noindent (b) $x$ has a continuous probability distribution for each
 $\theta$; 

\noindent (c) for each $x$ the probabilities in (b) are continuous
functions of $\theta$.

[We note that this definition includes sufficient statistics,
which are well known to be treatable as a single observation.]

\bigskip

Figure 1 illustrates application of the geometric formulation to 
determining the fiducial distribution (FD), $F_f(\theta|x)$, when the
random distributions (RD), $F_r(x|\theta)$, are non-intersecting.
For each fixed value of $x(=x_0$ say), the vertical
line intersects the RD's for each $\theta$ at a
probability value $F_r\equiv F_r(x_0|\theta)$, giving the points 
$(F_r,\theta)$. Non-intersection of the RD's implies,
as proven later (Theorem 1 in Subsection 3.2), 
that the $F_r$ value varies monotonically with $\theta$. 
The $(F_r,\theta)$ pairs, when plotted as $F_f(\equiv F_r)$ vs. $\theta$
(as shown in the right part of Figure 1) then constitute the fiducial
distribution, $F_f(\theta|x_0)$, provided that the RD's satisfy also the
following {\em completeness condition}: 
At the extreme values of $\theta$ the limiting RD's are 
identically zero and one, respectively, for all $x$ (Theorem 2 in
Subsection 3.2).  
 
\begin{figure}
 \psfig{file=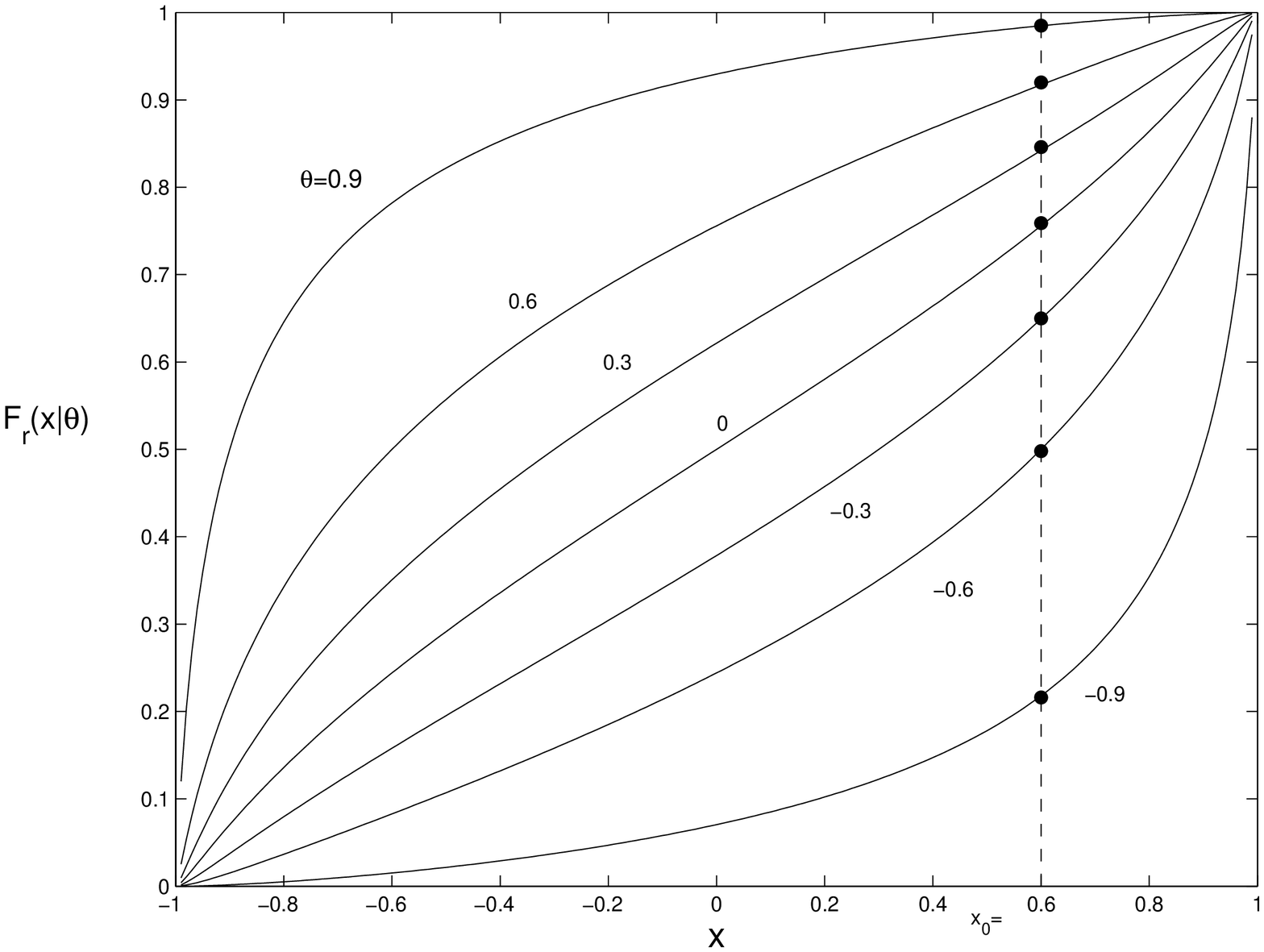,width=9.55cm,angle=0}
\hskip -0.9cm
\psfig{file=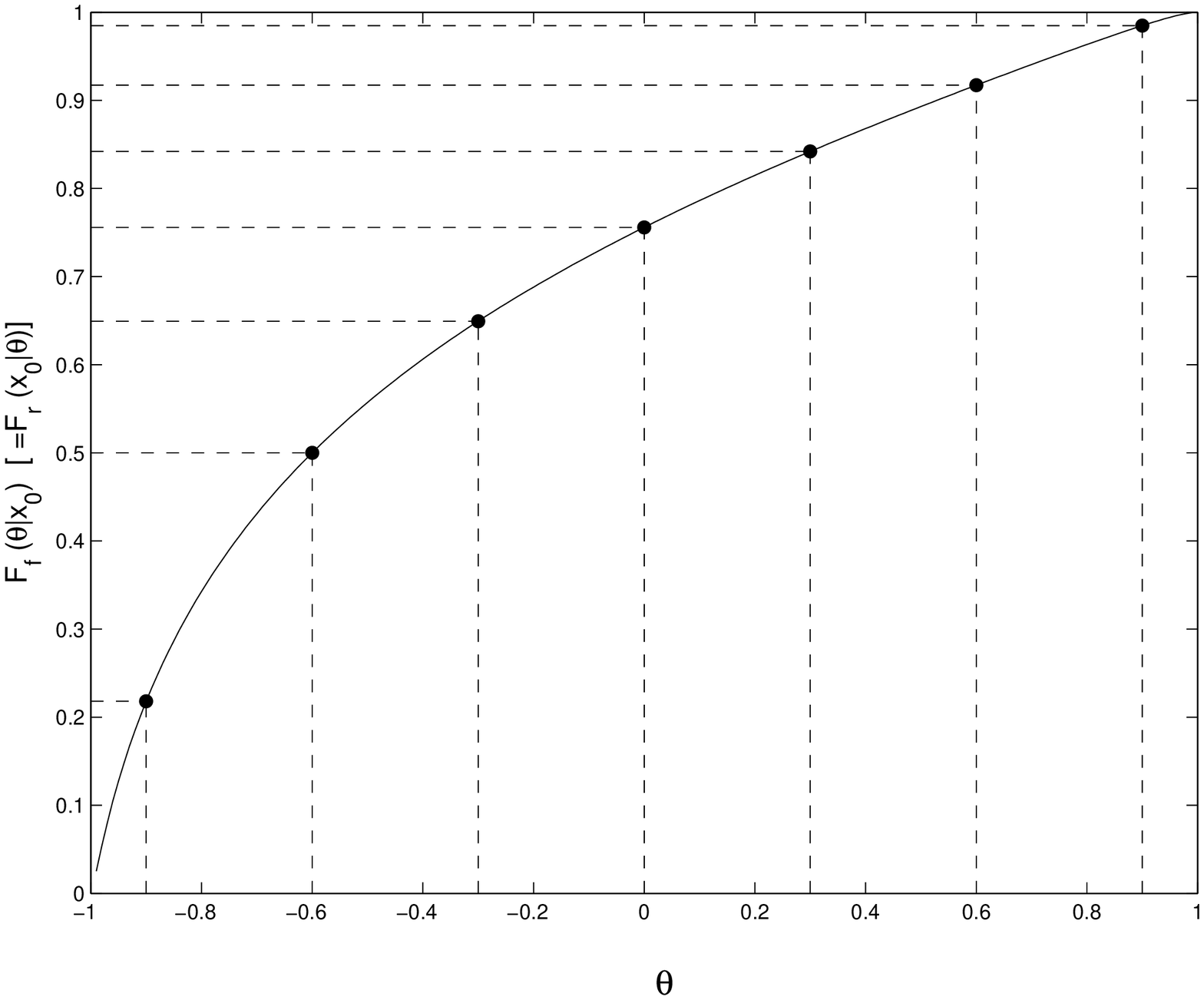,width=9.69cm,angle=0} 
\caption{FD $F_f(\theta|x_0)$ derived from RD $F_r(x_0|\theta)$}
\end{figure} 
  
This construction motivates the following definition, which is  justified
by the theorems in Section 3 cited above.

\noindent{\bf Definition CFD} {\em The (continuous) fiducial distribution
for each $x$,
$F_f(\theta|x)$, is the distribution induced on $\theta$ by the random
probability distributions (RD), $F_r(x|\theta)$, when (i) the RD's 
for different $\theta$'s are non-intersecting, and (ii) the RD's are
complete, i.e. approach zero and one, respectively, for all $x$ as
$\theta$ approaches its two extreme values.}

\noindent The continuous fiducial model then becomes an {\it FD model};
otherwise, if
either condition (i) or (ii) is not satisfied, we have a {\it non-FD
model}. 

Lindley's analysis in [7] of bayesian and fiducial
distributions assumes condition (ii) and the monotone conclusion of (i).

The above construction, which defines the FD $F_f(\theta|x_0)$ as
equivalent to the
pairs, $(F_r(x_0),\theta)$, also describes a {\em geometric identity}. 
For the situation in 
Figure 1 where $F_r$ is monotone increasing with $\theta$ (also referred to
as $\theta$ increasing), 
the equation  for this identity is 
\begin{equation}\label{eq1}
F_f(\theta|x_0)=F_r(x_0|\theta)\ \ (F_r \mbox{ increasing with } \theta).
\end{equation}
(In practice, this geometric identity is essentially
equivalent to the geometric formulation.)
At the minimum (infimum) value $\theta_m$ and maximum (supremum) value
$\theta_M$ the completeness condition is: For all $x$, 
$$\lim_{\theta\to\theta_m} F_r(x|\theta)=0,\\ 
\lim_{\theta\to\theta_M}F_r(x|\theta)=1.$$
In almost all applications, however, $F_r$ is monotone decreasing with
$\theta$.
($x_0$ then becomes the approximate median of the FD, as 
in the usual representation of translation parameters.) 
Replacing $F_f$ by $1-F_f$ yields 
\begin{equation}\label{eq2}
F_f(\theta|x_0)=1-F_r(x_0|\theta)\ \ (F_r \mbox{ decreasing with } \theta).
\end{equation}
(A corresponding change also occurs in the above completeness equations.)
The density
form of Eq. (\ref{eq2}) is also Fisher's
equivalence
equation [3, p.70]. (Fisher apparently did not consider using it to replace
the PQ formulation.)

Equations (\ref{eq1}) and (\ref{eq2}) are equivalent in the sense that
replacing 
$\theta$ by $-\theta$ (or by $1/\theta$ when the $\theta$ domain is
positive) converts either equation to the other. The symmetric
Eq.(\ref{eq1}) yields simpler analyses and is used throughout this paper, 
except for the example in Section \ref{section2}.

The identity Eq.(\ref{eq2}) can be expressed in terms of
fiducial
probability $P_f()$ and random probability $\Pr()$. With $X$ and $\Theta$ 
having the customary meaning associated with random variables, we have:
\begin{equation}\label{eq135}
Pr(X\geq x_0|\theta)=P_f(\Theta\leq \theta|x_0).
\end{equation}
[The term ``likelihood'', a common synonym for 
fiducial probability, probably represents the most appropriate 
interpretation to be associated with fiducial probability.]

\bigskip

\noindent{\bf Non-monotone measure for intersecting RD's}

The previous procedure for non-intersecting RD's in Figure 1 yielded the
points $(F_r,\theta)$ that comprised a monotone function. This function can
also be viewed as an ordinary monotone measure.  The relation is similar to
that when the RD's are complete: The monotone function becomes a fiducial
distribution $F_f(\theta|x_0)$ which then 
defines a probability measure of sets from generalizations of
the right side of Eq.(\ref{eq135}).

The same procedure is applicable to intersecting  RD's. The points
$(F_r,\theta)$ then comprise a non-monotone function, denoted by
$m_f(\theta|x_0)$, which represents also an associated {\it non-monotone
measure}. This is also referred to as a signed measure, being positive
(say) when $m_f$ is increasing and negative when decreasing. (A detailed
treatment of signed measures is given in Halmos [4, Chapter 6].)

In addition, the fiducial distribution notation, $F_f(\theta|x_0)$, in the
geometric identity Eq.(\ref{eq1}) is replaced by $m_f()$ giving a {\em
non-monotone geometric identity}:
\[m_f(\theta|x_0)=F_r(x_0|\theta).\] 

We also employ the notation $M_f()$ to denote either a monotone or a 
non-monotone
function, and with  an associated {\it general fiducial measure (FM)} that
can be either a monotone measure (when RD's are non-intersecting) or a 
non-monotone measure (when RD's are intersecting). ($M_f$ is especially
appropriate in the Section 3 proof prior to demonstrating existence of the
positive measure FD.) The {\em generalized geometric identity} then
becomes:
\begin{equation}\label{eq3}
M_f(\theta|x_0)=F_r(x_0|\theta).
\end{equation}

The same $(F_r,\theta)$ procedure also yields the following simple, but
important, proposition:

\noindent{\bf Proposition FM}\\
{\em For any family of parametric RD's, both the general fiducial measure
$M_f(\theta|x_0)$ and the generalized geometric identity Eq.(\ref{eq3})
always exist.}

A formal proof in terms of the surface $F(x,\theta)$ is given in Section 4.
In brief, existence of $M_f$ follows from the fact -- applicable to any
function $f(x,y)$ of two variables -- that the functions 
of one variable, say $f(x|y)$, determine the related surface, which then
also
determines  the functions of the other variable,
$f(y|x)$. The relationship
between the two functions is the geometric identity. When the values of
$f(x|y)$ are RD probabilities one gets measure terminology.

We note that item (c) in the continuous fiducial model Definition CFM can
now be
restated as: {\em (c') $M_f(\theta|x)$ is continuous for each $x.$} 
The monotone FM conclusion of Theorem 1 in Section 3 requires only this
continuity condition and the general geometric identity Eq.(\ref{eq3}). We
observe
also that Definition CFD becomes the statement: The general fiducial
measure $M_f(\theta|x_0)$
becomes a fiducial probability distribution $F_f(\theta|x_0)$ iff the RD's
satisfy conditions (i) and (ii).

\bigskip

\section{Non-FD joined distribution}\label{section2}

In connection with the proofs in Section \ref{section3}, it is useful to
first consider a non-FD example. A simple such example can be obtained by
joining two
different RD types. Assuming an unknown (decreasing) translation parameter
$\theta$, a normal
distribution for $\theta>0$ could be joined to (say) a Cauchy distribution
for $\theta <0$, or
more simply to merely another normal distribution with a different known
standard deviation. For computational ease we join two
families of uniform
distributions with different semiranges (equal to one-half the range) 
(see Figure 2(a)).
Continuity can be achieved by
introducing a $\theta$ {\it transition interval}, 
$[-\theta_T,+\theta_T]$, with transition RD's that vary monotonically
with $\theta$ from the RD at $+\theta_T$ to the RD at $-\theta_T$.
(Otherwise, as shown below, one gets a discontinuity at $\theta=0$.)
The semirange is $b$ for $\theta<-\theta_T$
and $a$ $(a<b)$ for $\theta>\theta_T$.

We assume a linear variation of the semiranges 
in the transition interval, $-\theta_T<\theta<+\theta_T$, yielding the
semiranges $T(\theta)$:
\begin{equation}\label{eq4}
T(\theta)=b[1-\theta/\theta_T]/2 +a[1+\theta/\theta_T]/2.
\end{equation}
The RD cumulative probability corresponding to a semirange $S$, which is
$a,b$ or $T(\theta)$ depending on the appropriate $\theta$
interval, is:
\[F_r(x|\theta)={1\over2}+\frac{x-\theta}{2S},\ \ \ -S\leq x-\theta\leq
S.\]

We note that all $\theta$ RD's in the transition interval pass through 
the point
$x_T$ where the $a$ and $b$
RD's for $\pm \theta_T$, respectively, intersect (see Figure 2(a)). We have
\[x_T=\frac{b+a}{b-a}\theta_T,\]
with corresponding probability 
\[F_T\equiv F_r(x_T|\pm\theta_T)={1\over2}+{\theta_T\over b-a}.\]
$(x_T, F_T)$ is also the vertex of a conical {\em intersection
region} with sides formed by these same $a$ and $b$ RD's.

\begin{figure}
\begin{center}
\setlength{\unitlength}{0.00059333in}
\begingroup\makeatletter\ifx\SetFigFont\undefined%
\gdef\SetFigFont#1#2#3#4#5{%
  \reset@font\fontsize{#1}{#2pt}%
  \fontfamily{#3}\fontseries{#4}\fontshape{#5}%
  \selectfont}%
\fi\endgroup%
{\renewcommand{\dashlinestretch}{30}
\begin{picture}(8579,6798)(0,-10)
\put(5712,5175){\blacken\ellipse{76}{76}}
\put(5712,5175){\ellipse{76}{76}}
\path(4212,450)(4212,6450)
\path(12,450)(8412,450)
\path(4212,450)(6612,6450)
\path(4799,450)(7199,6450)
\path(5397,450)(7797,6450)
\path(5997,450)(8397,6450)
\path(8412,2100)(7812,450)
\path(8412,3600)(7212,450)
\path(8412,5100)(6612,450)
\path(1212,450)(7212,6450)
\path(8412,6450)(12,1275)
\path(12,2025)(7212,6450)
\path(12,2750)(6087,6480)
\path(12,3150)(5412,6450)
\path(12,3450)(4812,6450)
\path(12,3825)(4212,6450)
\dashline{60.000}(4812,150)(4812,6450)
\path(3597,495)(6072,6435)
\path(2412,465)(6612,6420)
\dashline{60.000}(5712,150)(5712,6450)
\path(57,1675)(7782,6475)
\path(12,450)(7737,6450)
\path(12,2400)(6612,6450)
\put(4137,225){\makebox(0,0)[lb]{{\SetFigFont{12}{14.4}{\rmdefault}{\mddefault}{\updefault}{\tiny
0}}}}
\put(5412,225){\makebox(0,0)[lb]{{\SetFigFont{12}{14.4}{\rmdefault}{\mddefault}{\updefault}{\tiny1}}}}
\put(6612,225){\makebox(0,0)[lb]{{\SetFigFont{12}{14.4}{\rmdefault}{\mddefault}{\updefault}{\tiny2}}}}
\put(7812,225){\makebox(0,0)[lb]{{\SetFigFont{12}{14.4}{\rmdefault}{\mddefault}{\updefault}{\tiny3}}}}
\put(2937,225){\makebox(0,0)[lb]{{\SetFigFont{12}{14.4}{\rmdefault}{\mddefault}{\updefault}{\tiny-1}}}}
\put(1662,225){\makebox(0,0)[lb]{{\SetFigFont{12}{14.4}{\rmdefault}{\mddefault}{\updefault}{\tiny-2}}}}
\put(537,225){\makebox(0,0)[lb]{{\SetFigFont{12}{14.4}{\rmdefault}{\mddefault}{\updefault}{\tiny-3}}}}
\put(3850,6375){\makebox(0,0)[lb]{{\SetFigFont{12}{14.4}{\rmdefault}{\mddefault}{\updefault}{\tiny1.0}}}}
\put(4002,975){\makebox(0,0)[lb]{{\SetFigFont{12}{14.4}{\rmdefault}{\mddefault}{\updefault}{\tiny.1}}}}
\put(4002,1575){\makebox(0,0)[lb]{{\SetFigFont{12}{14.4}{\rmdefault}{\mddefault}{\updefault}{\tiny.2}}}}
\put(4002,2175){\makebox(0,0)[lb]{{\SetFigFont{12}{14.4}{\rmdefault}{\mddefault}{\updefault}{\tiny.3}}}}
\put(4002,2775){\makebox(0,0)[lb]{{\SetFigFont{12}{14.4}{\rmdefault}{\mddefault}{\updefault}{\tiny.4}}}}
\put(4002,3450){\makebox(0,0)[lb]{{\SetFigFont{12}{14.4}{\rmdefault}{\mddefault}{\updefault}{\tiny.5}}}}
\put(4002,3975){\makebox(0,0)[lb]{{\SetFigFont{12}{14.4}{\rmdefault}{\mddefault}{\updefault}{\tiny.6}}}}
\put(4002,4575){\makebox(0,0)[lb]{{\SetFigFont{12}{14.4}{\rmdefault}{\mddefault}{\updefault}{\tiny.7}}}}
\put(4002,5175){\makebox(0,0)[lb]{{\SetFigFont{12}{14.4}{\rmdefault}{\mddefault}{\updefault}{\tiny.8}}}}
\put(4002,5775){\makebox(0,0)[lb]{{\SetFigFont{12}{14.4}{\rmdefault}{\mddefault}{\updefault}{\tiny.9}}}}
\put(3262,1000){\makebox(0,0)[lb]{{\SetFigFont{12}{14.4}{\rmdefault}{\mddefault}{\updefault}{\tiny$\theta_T$=.5}}}}
\put(7750,5950){\makebox(0,0)[lb]{{\SetFigFont{12}{14.4}{\rmdefault}{\mddefault}{\updefault}{\tiny$\leftarrow\theta_T$=-.5}}}}
\put(5370,6200){\makebox(0,0)[lb]{{\SetFigFont{12}{14.4}{\rmdefault}{\mddefault}{\updefault}{\tiny$\theta_T$=.5$\to$}}}}
\put(4687,1950){\makebox(0,0)[lb]{{\SetFigFont{12}{14.4}{\rmdefault}{\mddefault}{\updefault}{\tiny1}}}}
\put(5887,1950){\makebox(0,0)[lb]{{\SetFigFont{12}{14.4}{\rmdefault}{\mddefault}{\updefault}{\tiny2}}}}
\put(7087,2000){\makebox(0,0)[lb]{{\SetFigFont{12}{14.4}{\rmdefault}{\mddefault}{\updefault}{\tiny3}}}}
\put(8000,2000){\makebox(0,0)[lb]{{\SetFigFont{12}{14.4}{\rmdefault}{\mddefault}{\updefault}{\tiny$\theta$=4}}}}
\put(5407,0){\makebox(0,0)[lb]{{\SetFigFont{11}{13.2}{\rmdefault}{\mddefault}{\updefault}{\tiny$x=1.25$}}}}
\put(4462,0){\makebox(0,0)[lb]{{\SetFigFont{11}{13.2}{\rmdefault}{\mddefault}{\updefault}{\tiny$x=0.5$}}}}
\put(8487,450){\makebox(0,0)[lb]{{\SetFigFont{11}{13.2}{\rmdefault}{\mddefault}{\updefault}$x$}}}
\put(3912,6675){\makebox(0,0)[lb]{{\SetFigFont{11}{13.2}{\rmdefault}{\mddefault}{\updefault}$F_r(x|\theta)$}}}
\put(2707,1305){\makebox(0,0)[lb]{{\SetFigFont{12}{14.4}{\rmdefault}{\mddefault}{\updefault}{\tiny.25}}}}
\put(5577,5210){\makebox(0,0)[lb]{{\SetFigFont{11}{13.2}{\rmdefault}{\mddefault}{\updefault}{\tiny $A$}}}}
\put(5677,5140){\makebox(0,0)[lb]{{\SetFigFont{11}{13.2}{\rmdefault}{\mddefault}{\updefault}$\bullet$}}}
\put(5030,4520){\makebox(0,0)[lb]{{\SetFigFont{11}{13.2}{\rmdefault}{\mddefault}{\updefault}{\tiny $x_T$}}}}
\put(5197,4450){\makebox(0,0)[lb]{{\SetFigFont{11}{13.2}{\rmdefault}{\mddefault}{\updefault}$\bullet$}}}
\put(2667,2095){\makebox(0,0)[lb]{{\SetFigFont{12}{14.4}{\rmdefault}{\mddefault}{\updefault}{\tiny0}}}}
\put(2520,2675){\makebox(0,0)[lb]{{\SetFigFont{12}{14.4}{\rmdefault}{\mddefault}{\updefault}{\tiny-.25}}}}
\put(2482,2985){\makebox(0,0)[lb]{{\SetFigFont{12}{14.4}{\rmdefault}{\mddefault}{\updefault}{\tiny$\theta_T$=-.5}}}}
\put(2592,3405){\makebox(0,0)[lb]{{\SetFigFont{12}{14.4}{\rmdefault}{\mddefault}{\updefault}{\tiny-1}}}}
\put(2577,4080){\makebox(0,0)[lb]{{\SetFigFont{12}{14.4}{\rmdefault}{\mddefault}{\updefault}{\tiny-2}}}}
\put(2547,4815){\makebox(0,0)[lb]{{\SetFigFont{12}{14.4}{\rmdefault}{\mddefault}{\updefault}{\tiny-3}}}}
\put(2300,5535){\makebox(0,0)[lb]{{\SetFigFont{12}{14.4}{\rmdefault}{\mddefault}{\updefault}{\tiny$\theta$=-4}}}}
\put(660,5535){\makebox(0,0)[lb]{{\SetFigFont{12}{14.4}{\rmdefault}{\mddefault}{\updefault}{\Huge{\bf(a)}}}}}
\end{picture}
}
\psfig{file=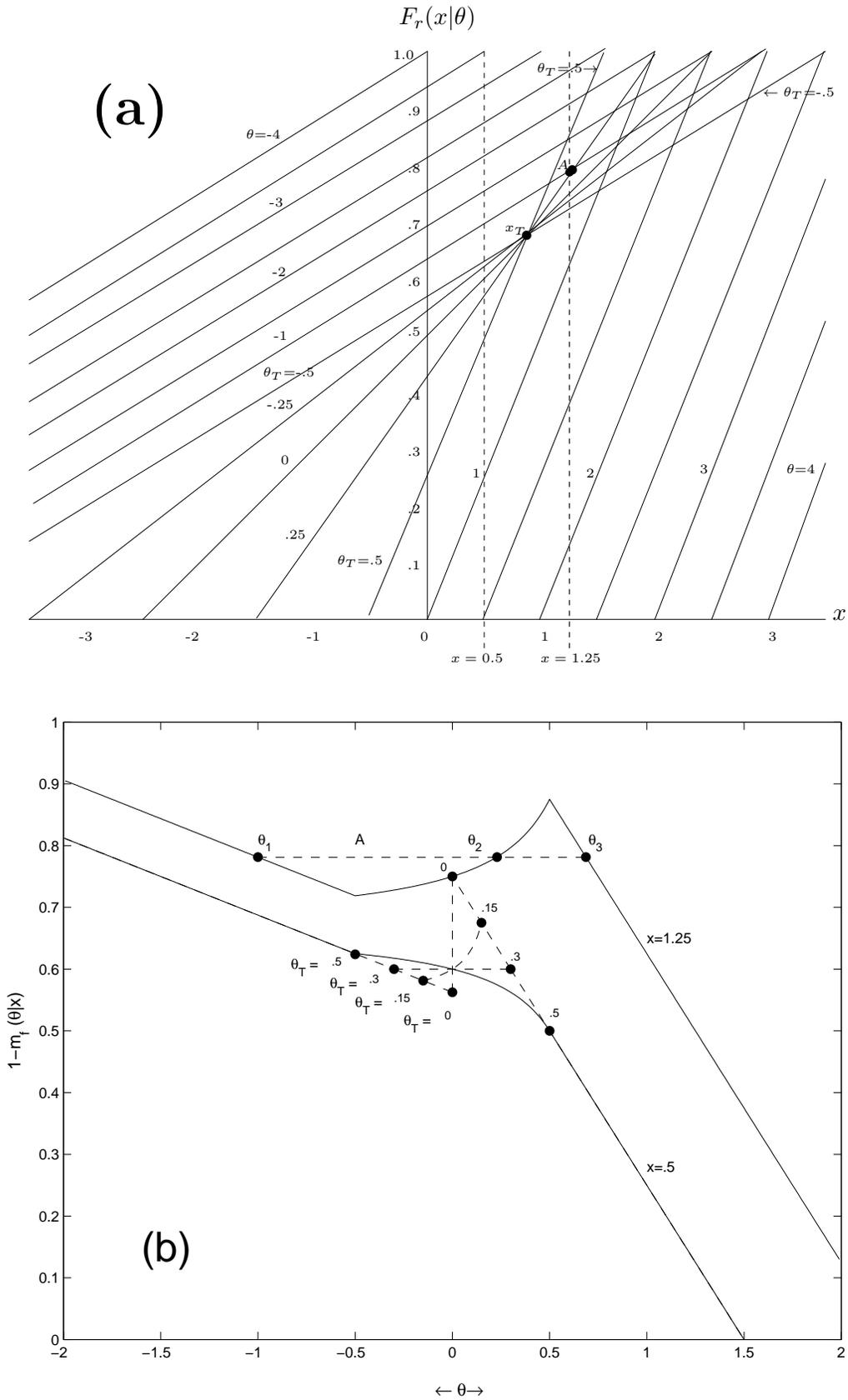,width=16cm,angle=0}
\vskip-0.7cm
\caption{(a) Non-FD joined distributions, (b) non-monotone FM for
intersecting RD's} \label{f2}
\end{center}\end{figure}

Figure 2(a) plots the RD's for numerical values $a=1$, $b=4$,
$\theta_T$=.5,
yielding the transition intersection point $x_T=5/6$ with probability
$F_T$=2/3. The RD intersection region is defined by $x>(5/6)$ and 
$\theta \in [-.5,+.5]$. At each point $(x,F_r)$ in this region, three
different
RD's intersect, corresponding to semirange $a$, semirange $b$ and
transition semirange $T(\theta)$.

The intersections are shown more clearly in Figure 2(b), which plots
the (signed) fiducial measure for $x=1.25~(>x_T)$ and $ x=.5~(<x_T)$. The
ordinate is $1-m_f(\theta|x)$ rather than $m_f(\theta|x)$ in
order to show more directly the relationship with the RD's in
Figure 2(a); since $\theta$ is a decreasing parameter, it is now 
$1-m_f(\theta|x)$
that equals $F_r$. 
The 
fiducial curve for $x=1.25$
consists of three parts: The monotone decreasing linear part on the left 
corresponds in Figure 2(a) to the values of $F_r$ as one proceeds 
downward on the $x=1.25$ vertical line with semirange $b$ RD's from 
$\theta=-2.75$ to $\theta=-.5(=-\theta_T)$, the lower boundary of the 
intersection region. For the second intermediate part, one has a monotone
increasing transition curve with variable transition semiranges from 
$\theta =-\theta_T$ to $+\theta_T$; these correspond in
Figure 2(a) to the $F_r$ values as one proceeds upward from the lower 
($-\theta_T$) to the upper ($+\theta_T$) transition boundary. The third
monotone decreasing linear part, similar to the first part, corresponds in
Figure 2(a) to the values of $F_r$ as one proceeds downward with
semirange $a$ RD's from $\theta=.5(=+\theta_T)$, the upper boundary of 
the intersection region, to $\theta=2.25$.

In Figure 2(a) the point (labeled $A$) at
$x=1.25$ and $ F_r=0.78$, corresponds on the Figure 2(b) $x=1.25$
curve to fiducial measure $1-m_f=0.78$ (the horizontal line $A$), and
yields the three intersections:  $\theta_1=-1$ for the RD with 
semirange $b=4;\; \theta_2=.25$ for the
transition RD with semirange $T(.25)=1.75$; and $\theta_3=.68$ for the RD
with semirange $a=1$. 

The 1-$m_f$ fiducial curve for $x=.5$, being less than the intersection
$x_T=.83$, is monotone decreasing for all $\theta$. 
Also shown (by dashed lines) is the effect of reducing the transition
interval $[-\theta_T,+\theta_T]$. When $\theta_T=.3$ the intersection point
$x_T=.5$ is the same value used to obtain the $x=.5$ fiducial curve. Hence
the fiducial transition 
is constant, reflecting the fact that the vertical line in
Figure 2(a) passes through $x_T$. Further reduction to the $\theta_T =.15$
curve 
places the $x=.5$ line of Figure 2(a) in the intersection region with an
increasing fiducial transition curve. As $\theta_T$ and $x_T$
approach zero the fiducial transition
degenerates into a discontinuous vertical jump. The final paragraphs below
discuss further the discontinuous solution.

We emphasize several features of this example. Depending on
the value of $x$, there exist both
non-intersecting RD regions where the fiducial measure $m_f(\theta|x)$ 
is monotone and hence constitutes a formal distribution, and also 
intersection regions where $m_f$ is non-monotone and hence not a
distribution. (The non-FD label affirms that a fiducial 
distribution analysis is
meaningful only if applicable to all admissible observations.)
 Another feature is the existence, for values of 
$x$ where intersections occur, of multiple
values of $\theta$ having the same fiducial measure $m_f(\theta|x)$. 
This is characteristic
of a non-monotone measure, and is in fact a necessary and sufficient
condition for non-monotonicity (Theorem EQ/NM in Appendix B).

\bigskip

\noindent{\bf The discontinuous solution}

The transition interval has degenerated to a point
$\theta_T=0$ and hence contains no RD. However, the intersection
region still exists, bounded by the semirange $a$ and $b$ RD's and by
$x>0$.
In this region each value of $m_f$ still yields three intersections, 
corresponding to semirange $a$ and $b$ RD's plus the intersection with the 
vertical jump. Although this last intersection is not associated with any
RD, one can obtain a {\it normalized transition RD$^*$} corresponding to
the normalized parameter
$\theta^*\equiv \theta/\theta_T,$ $(-1\leq \theta^*\leq +1)$. In fact,
the previous Eq.(\ref{eq4}) for semirange $T(\theta)$ can be 
interpreted as
a formula for a normalized semirange $T^*(\theta^*)$ and corresponding
normalized RD$^*$. Together they constitute a {\it transition
protocol} that is applicable to all $\theta_T$. The actual non-normalized
transition RD's and associated
$\theta$'s for a specific non-zero $\theta_T$ are obtained from
$\theta=\theta^*\cdot\theta_T.$ In the limit the normalized
$\theta^*$'s apply also to $\theta_T=0$, and a normalized RD$^*$ can then
be associated with the $m_f$ intersection on the vertical discontinuity.

This RD$^*$ solution for $\theta_T=0$ can be viewed in either a negative
or a positive light. For the former the RD$^*$ value obtained  depends upon
an assumed transition protocol, although the discontinuous formulation with
joining at $\theta=0$ does not involve any transition. Consequently
RD$^*$ is essentially indeterminate. The opposing positive view is that
associating
the $\theta=0$ discontinuity with the limit of {\it any} transition
protocol illuminates the underlying nature of the discontinuity. It is not
really necessary to assign a specific RD$^*$ value. 

\section{Fiducial non-intersection theorems}\label{section3}

This section proves  necessary and sufficient conditions for an FD
 to exist. These are implied assumptions in almost all applications.
\bigskip

\noindent{\bf Theorem 1} {\em
The fiducial measures (FM) $M_f(\theta|x)$ are monotone for
each $x$ iff the RD's $F_r(x|\theta)$ are non-intersecting for all
$\theta$.}

\noindent(For brevity we often write $M(\theta)\equiv M_f(\theta|x)$.)

\bigskip

\subsection{Strict non-intersection}
We consider first the simpler case of {\it strictly non-intersecting}
RD's, i.e. completely separated with different $\theta$ RD's having no
points in common, together with corresponding {\it strictly monotone}
FM's, i.e. no  intervals where FM is constant. (The complete non-strict
case is treated below in Theorem 1.)

Strict non-intersection does not
apply to finite endpoints $x_m$ and $x_M$ where the same probability,
 $F_r=0$ or 1, applies to all RD's and the corresponding FM
$M_f(\theta|x_{m/M})$ is constant. More generally, this {\it endpoint
condition} exists when an intersection point is common to all $\theta$.

\bigskip

\noindent{\bf Theorem 1S} {\em   
$M_f(\theta|x)$ is strictly monotone for each $x$ iff the RD's 
$F_r(x|\theta)$ for different $\theta$ are strictly non-intersecting.}
\bigskip

We recall, and cite for reference, the following well known definition:
{\em $M(\theta)$ is a strictly monotone increasing (say) function
iff for all $\theta_2 > \theta_1$ we have $M(\theta_2)> M(\theta_1)$.}


This definition is of course a special case of the general monotone
definition wherein the above $>$ relation 
is replaced by $\geq$, giving $M(\theta_2)\geq M(\theta_1)$.
Note that these definitions allow discontinuous jumps. It is more
convenient, however, to introduce continuity when applying 
Lemma 1 below.

The lemma is a generalization of the intermediate value
property of a continuous function.
For example, for any non-monotone function we have: if points 1 and 3 are
on opposite sides and in a neighborhood of a (point 2) maximum (say),
then one can draw a horizontal line between
point 2 and points 1, 3.
\bigskip

\noindent{\bf Lemma 1} {\em
Let $f(z)$ be a continuous function. Suppose the ordered points
$z_1<z_2<z_3$ have corresponding values of $f(z_i)$ that are not in
monotone order, i.e. $f(z_1)<f(z_2)>f(z_3)$ or $f(z_1)>f(z_2)<f(z_3)$. 
Then there exist points $z_1^*, z_2^*$  such that $f(z_1^*)= f(z_2^*).$}

\begin{proof} For the case $f(z_1)<f(z_2)>f(z_3)$ let $f^*$ be such
that $f(z_2)>f^*>\max\{f(z_1), f(z_3)\}$. For the continuous segment
$(z_1,z_2)$ the intermediate value theorem for a continuous function yields
$z_1^*$ such that $f(z_1^*)=f^*$. Likewise, for the continuous segment
$(z_2,z_3)$ there exists $z_2^* (>z_1^*)$ such that $f(z_2^*)=f^*$. The
proof for the other non-monotone case is similar.
\end{proof}
\bigskip

For later application we use the following simple extension of Lemma 1:

\bigskip

\noindent{\bf Corollary} {\em 
In addition to the assumptions of the lemma, let $M_1,M_2,\ldots$ be 
any finite or denumerable set of numbers.
Then ${z'}_1^{*}, {z'}_2^{*}$ exist such that both
$f({z'}_1^{*})=f({z'}_2^{*})$, and which are also  not equal to $M_i\;
(i=1,2,\ldots)$.}

 \noindent\begin{proof}
The set of admissible function values $f^*$ in the proof of the lemma
includes an open set $\cal F$, say. The difference subset 
${\cal F}\setminus\{M_1,M_2,\ldots\}$ subtracts at most a denumerable set
from $\cal F$ whence the remainder is non-empty.
\end{proof}

\bigskip

\begin{proof11}   Suppose that the continuous RD's $F_r(x|\theta)$ are
strictly non-intersecting. We initially presume, until monotonicity 
of the FM with
respect to $\theta$ is proven, that a $\theta$ numerical value (in the
interval domain specified as assumption (a) in the Section 1 fiducial model
Definition CFD) is
attached to each different RD and constitutes an
arbitrary ``label'' for the RD. Then for each $x$ and
all $\theta_i\neq \theta_j$, strict non-intersection implies that
$F_r(x|\theta_i)\neq F_r(x|\theta_j).$ Equivalently, in terms of the FM
we have, from the geometric identity Eq.(\ref{eq3}), that $M(\theta_i) 
\neq M(\theta_j)$. We need to prove that continuity of
$M(\theta)$ implies that $M(\theta)$ is strictly monotone. 

Let $\theta'$ and $\theta''(>\theta')$ be an arbitrary ``reference"
pair and suppose that $M(\theta'')>M(\theta')$. (The opposite case,
$M(\theta'') <M(\theta')$, is treated similarly.) For any
$\theta_1,\theta_2$ pair with $\theta_2 >\theta_1$ (not necessarily
different from $\theta'$ or $\theta''$) we will show that
continuity of $M(\theta)$ implies that $M(\theta_2)>M(\theta_1)$. Suppose
the contrary, that $M(\theta_2) <M(\theta_1)$, and consider the six
different inequality cases satisfied by 
the four  $\theta$'s,  in combination with the various corresponding
  relative positions (inequalities) of the four quadruplet values,
$M(\theta_1), M(\theta_2),
M(\theta'), M(\theta'')$ (abbreviated to $M_1,M_2,M',M''$), that satisfy
the two assumed inequalities. It is readily
seen that at most three of the four M's can be monotonically related;
moreover, an appropriate selection of two $M$'s from any such monotone
triplet, together with the
fourth $M$, will comprise a non-monotone $M$ triplet. 

For example, consider the case $\theta'<\theta_1< \theta_2 <\theta''.$ If
$M_1<M'$, then the triplet $M'> M_1>M_2$ is monotone decreasing and
selecting any two of these
together with $M''$, is a non-monotone triplet. (If $\theta_1=\theta'$  the
resultant triplet is already non-monotone.) For
the case where $M_1> M''$ and $M_2<M'$, all triplets are
non-monotone. 

Thus the requirements of Lemma 1 are satisfied, and there exist
$\theta^*_1\neq\theta^*_2$ with $M(\theta^*_1)= M(\theta^*_2)$,
which contradicts the non-intersection condition that $M(\theta^*_1) \neq
M(\theta^*_2).$ Hence $M(\theta_2)>M(\theta_1)$ and $ M(\theta)$
is strictly monotone increasing. Similarly, if initially
$M(\theta'')<M(\theta')$, then $M(\theta)$ is strictly monotone decreasing.

For the converse, $M_f(\theta|x)$ is assumed for each $x$ to be
either strictly monotone increasing or decreasing. If the former,
then for all $x,\;\; \theta_2>\theta_1$ implies that 
$M(\theta_2)>M(\theta_1)$. That is, the RD's at each $x,\; F_r(x|\theta)$, 
are separated for different $\theta$'s.  Combining these $F_r(x|\theta)$'s
for varying $x$ yields separated
and hence strictly non-intersecting RD's. Similarly for monotone
decreasing.  
\end{proof11}

\bigskip

The four $\theta$'s and inequalities in the proof are also the same as
in Definition NM in Appendix B of an (assumed contrary) non-monotone
$M(\theta)$.  We also remark that a variation of the proof merely cites the
non-monotone Theorem EQ/NM, also in Appendix B, wherein a continuous
non-monotone $M(\theta)$ implies existence of $M(\theta_1)=M(\theta_2)$,
and
thereby providing the contradiction needed to prove Theorem 1S. However,
the proof of Theorem EQ/NM is essentially the same as for Theorem 1S and
also requires the use of Lemma 1.

We note also that the assumption of FM
continuity (item (c) or (c') in Definition CFM) can be relaxed to permit a 
denumerable number of jumps. This would be applicable, however, only for
RD's with
both a continuous and a discrete component.

\subsection{General non-intersection}

\noindent{\bf Touching}

Different RD's may be non-intersecting even though they have
points in common. This occurs when RD's {\it touch} at a
point $x_0$, say.
For two RD's with $\theta_1<\theta_2$ we have: (i) $F_r(x_0|\theta_1)=
F_r(x_0|\theta_2)$, and (ii) for right and left
(non-touching) neighborhoods of
$x_0$, with elements denoted by ${x+}$ and ${x-}$, respectively, the
monotone directions with $\theta$ are the same, i.e. for monotone
increasing (say) with $\theta$, we have both $F_r({x+}|\theta_1)
<F_r({x+}|\theta_2)$ and $F_r({x-}|\theta_1) <F_r({x-}|\theta_2)$. 

\noindent(When RD derivatives exist, the RD's are tangent and unequal 
convexities (say)
yield (ii).) (An intersection satisfies condition (i), but in (ii) the
monotone directions in right and 
left neighborhoods are opposite.)

Touching at $x_0$ necessarily includes all
intermediate $\theta$ in the closed interval 
$[\theta_1,\theta_2].$ 
The corresponding RD's are also, because  of
continuity, intermediate to the $\theta_1$ and $\theta_2$ RD's in
the right and left neighborhood inequalities, and hence necessarily
pass through $x_0$. (A formal proof is included in the proof below of
Theorem 1.) (A different route is taken by intersecting RD's in
reaching the reverse inequality in (ii); see the Section \ref{section2}
example.)

$[\theta_1,\theta_2]$ is a maximum size interval (assumed hereafter)
iff all $\theta$ that are either $<\theta_1$ or $>\theta_2$ are
non-touching in the neighborhoods of $x_0.$
$M_f(\theta|x_0)$ is constant for $\theta\in[\theta_1,\theta_2].$
Indeed, existence of touching intervals is
necessary and sufficient for a FM to be non-strict, and also for
non-intersecting RD's to be non-strict.

For the same $x_0,\; M_f(\theta|x_0)$ may contain several touching
intervals. The $\theta$ intervals are disjoint since the $F_r$ values
are different. Between  two successive constant FM intervals is an open
interval of non-touching $\theta$'s with strictly 
monotone $M_f(\theta|x_0)$.

Additional touching properties are presented following the proof of Theorem
2.

\bigskip

\noindent
{\bf Proof of Theorem 1} 

We first restate the usual definition of general
monotonicity to allow explicitly for constant intervals, by separating
the $\geq$ relation into $=$ and $>.$
\bigskip
   
  \noindent{\bf Definition GM} {\em
$M(\theta)$ is monotone increasing (say) iff:
(i)  $M(\theta)=M_i, (i=1,\ldots)$ for $\theta$ in disjoint intervals
$I_i$; otherwise (ii) for any pair
$(\theta_1,\theta_2)$ not in the same $I_i$ {\rm (the set of such pairs
is denoted by ${\cal P}$)}, we have that $\theta_2>\theta_1$ implies 
$M(\theta_2)>M(\theta_1)$.}  

\noindent
As before, jumps are not precluded, in either (i) or (ii), but are later
eliminated when continuity of $M(\theta)$ is introduced.

It is convenient to relate (using notation (i') and (ii')) each part
of the proof to each part ((i) and (ii)), respectively, in the above
Definition GM. 
The non-intersection RD's $F_r(x|\theta)$, with arbitrary $\theta$
labeling (as in the proof of Theorem 1S), 
have for each $x,\; \theta$'s that are either: (i') within disjoint 
touching sets $E_i\; (i=1,\ldots,n_x)$, or are (ii') non-touching.
The non-intersection condition implies: (i') for the touching sets
$E_i$ we have $M(\theta)=M_i$ for all $\theta$ in $E_i$; and (ii') for all
pairs $\theta_1,\theta_2$ in ${\cal P}$, we have that
$\theta_2\neq\theta_1$ implies $M(\theta_2)\neq M(\theta_1)$. 

Re (i'): We need to show that continuity of $M(\theta)$ implies
that $E_1$ (say) is an interval $I_1$. If, to the contrary, $E_1$
were not an interval, there would exist one or more gaps. That is, there
exist $\hat \theta$ not in $E_1$ and a pair $(\theta_1, \theta_2)$ in $E_1$
such that $\theta_1<\hat\theta<\theta_2$. 
Now $M(\hat \theta)(\equiv \hat M) \neq M_1$, say $>M_1$ (the other $<$
case is similar), whence 
$M_1=M(\theta_1)< M(\hat \theta)>M(\theta_2)=M_1$. 
Since $M(\theta)$ is continuous, it follows from the Corollary to
Lemma 1 that for $M'^*$ satisfying both $M_1<M'^*<\hat M$ and also 
$M'^* \neq M_i,\;(i=2,\ldots)$, there will exist
${\theta'}^*_1$ and $ {\theta'}_2^{*}$ with
$M({\theta'}_1^*)=M({\theta'}_2^*)={M'}^*$.
Since the pair $({\theta'}^*_1, {\theta'}^*_2)$  is in $\cal P$
this contradicts the non-intersection
assumption that $M({\theta'}^{*}_1)\neq M({\theta'}^{*}_2)$ for all
$\theta$ pairs in ${\cal P}$. Hence $\hat \theta$ is in $E_1$, there 
is no gap, and $E_1$ is
an interval $I_1$. Similarly, for $M(\hat\theta)<M_1$ and for $i>1$.
Hence (i') becomes (i). 

Re (ii'):  The proof is practically the same as in Theorem 1S.
Let $\theta'$ and $\theta''(>\theta')$ be a reference pair in 
${\cal P}$ with (say) $M(\theta')<M(\theta''),$  and let 
$\theta_1$ and $ \theta_2 (>\theta_1)$ be any pair in ${\cal P}$.
Suppose, contrariwise, that $M(\theta_1)>M(\theta_2)$. For any
relative position of the four $\theta$'s and corresponding 
$M(\theta)$'s, then, as in the proof of Theorem 1S, 
there always exists a non-monotone $M(\theta)$ triplet. The 
Corollary to Lemma 1 then yields two ${\theta'}^*$'s in ${\cal P}$
with equal $M({\theta'}^*)$'s. Hence a contradiction and
$M(\theta_1)>M(\theta_2)$ in ${\cal P}$. Similarly for decreasing
monotonicity, $M(\theta') >M(\theta'')$. Hence $M(\theta)$ is monotone in
${\cal P}$, (i.e., for non-touching $\theta$'s) and (ii') becomes (ii).
Hence $M_f(\theta|x)$ is monotone for all $x$.

Conversely, let $M_f(\theta|x_0)$ be monotone increasing (say)
for each $x_0$. (a) When $M_f(\theta|x_0)=$constant for 
disjoint intervals
$[\theta_{1j}(x_0), \theta_{2j}(x_0)],$ $(j=1,\ldots,m_j(x_0))$,
then equivalently the RD's $F_r(x_0|\theta)$ are touching for
$\theta$'s in these touching intervals. (b) When $M_f$ is strictly
increasing in an interval, the equivalent $F_r(x_0|\theta)$ is also
strictly increasing with $\theta$. Combining both (a) and (b) $F_r$'s for
all $x$'s (replacing the $x_0$ notation) yields a non-intersecting
configuration. Similarly for $M_f$ monotone
decreasing, which proves the converse.
\endpf

There exists also an {\it endpoint solution} -- discussed in the
intersection
analysis in Appendix B -- corresponding to an endpoint
condition where all RD's intersect at the same point $x_I\neq 0,1$, 
and thereby constituting a common endpoint for two adjacent FM's. (An
example is a scale parameter for an asymmetrical distribution having both
positive and negative values.)
Although such a solution formally satisfies the monotone conditions for
Theorem 1, it is not applicable to the (implicit) single fiducial 
model and can be disregarded. 

We note the following symmetric relation between $F_r$ and $M_f$:

\noindent
$M_f(\theta|x)$ is [strictly] monotone for each $x$ iff
$F_r(x|\theta)$ is [strictly] non-intersecting for different $\theta$'s.
Contrapositively, $F_r(x|\theta)$ is
[strictly] monotone for each $\theta$ iff $M_f(\theta|x)$ is [strictly]
non-intersecting for each $x$.

Similar relations exist regarding touching and constant intervals; and also
completeness (in Theorem 2 below).
\bigskip

\noindent {\bf FD existence theorem}
\bigskip

  \noindent{\bf Theorem 2} {\em
A FD exists iff the RD's are non-intersecting and complete.}\\
\begin{proof}
  From Theorem 1, RD's are non-intersecting iff the FM's are
monotone for all $x$. Also, by definition, monotone FM's are (fiducial)
distributions iff the extreme values of $M_f(\theta|x)$ are 0 for
$\theta_m$ (say) and 1 for $\theta_M$, for all $x$, i.e. the RD's are
complete. (More accurately, for increasing (say) $\theta, \;
\lim_{\theta\rightarrow\theta_m} F_r(x|\theta)(=M_f(\theta|x))=0$ for all
$x$, and similarly $=1$ for $\theta_M$. 
These limiting distributions are
discontinuous at finite endpoints $x_M$ and $x_m$.) 
\end{proof}

Note that the endpoint solution does not satisfy the 0,1 limits.

Ordinarily, if non-intersecting RD's are incomplete, completeness can be
achieved by supplementing with additional non-intersecting RD's; the RD's
can then be said to be {\it completable}. An important exception is the
incompletable non-central chi-square distribution. (See Subsection 5.2)
In contrast to the non-intersection condition, completeness 
is a boundary-type condition. It is indicated in Subsection 5.2 that 
incompletable distributions occur only for certain non-FD 
``composite'' distributions.

\bigskip

\noindent{\bf Touching segments}

In addition to a touching $\theta$ interval at $x_0$, suppose there is also
touching at $x'_0 (>x_0)$. This  can occur
following an open right non-touching neighborhood about $x_0$,
which serves also as a left neighborhood about $x'_0$ for the new  
touching $\theta$'s.

Such touching can be continuously extended to an entire $x$-interval, 
say $[x_1, X]$. The $\theta$ interval
need not be the same for all $x$ in $[x_1, X]$.  Changes consist in
adding or removing $\theta$'s  while incorporating transition 
neighborhoods similar to those in successive point
touchings. The result is that changes can occur only at discrete points
$x_i$ $(i=1,\ldots,n)$, with an open $x$-interval between changes serving
as a transition. When adding $\theta$'s at $x_i$ the new
touching $\theta$
interval also begins at $x_i$, but when removing $\theta$'s the actual
reverting to non-touching status takes effect at $x_i+$, the ``beginning''
of the open right neighborhood of the removed $\theta$'s. [The
notation $x_i+$  (which differs slightly from the
previous usage) to denote a ``beginning'', is merely suggestive.
A rigorous condition states that $\theta$ is non-touching in all right
(sufficiently small) $\epsilon$ neighborhoods of $x_i$.]
When both addition and removal occur, only the added $\theta$'s are
included at $x_i$, 
while the touching $\theta$'s that remain after the removals appear at
$x_i+$. 
Changes in $\theta$'s can be to the beginning and/or 
the ending of the previous $[\theta_1, \theta_2]$ interval.
 
The beginning $x_1$ of an entire $x$-interval constitutes the initial
addition change of RD's from non-touching to touching status, while $x_n
(\equiv X)$ is the last touching before removal of all remaining RD's to
non-touching status at $x_{n}+$. Point touching at $x_0$ is the special
degenerate interval case, $x_1=x_n(=x_0)$, with $\theta$ touching
starting at $x_0$ and removal at $x_{0}+$.

The above changes have been viewed relative to increasing values of
$x$. The same result is obtained if $x$ decreases from $x_n$ to $x_1$;
previous additions are treated now as removals and previous removals
treated as additions. 

These results can be expressed as a {\it touching segment}  that consists 
of a {\it $\theta$ touching interval function}
$[\theta_L(x), \theta_U(x)]$ for $x$ in the segment domain
$[x_1,X]$. $\theta_L(x)$ and $\theta_U(x)$ are each step functions
with jump discontinuities at discrete points $x_1,\ldots,x_n(=X).$ 
At an $x_i$ discontinuity $\theta_U(x_i)$ assumes the larger of the two
values that are before and after the discontinuity. For adding $\theta$'s
this is the value after the positive jump; for removing $\theta$'s it is
before the negative jump. Both adding and removing (different) $\theta$'s
can be viewed as combining a positive jump at $x_i$ followed (almost)
simultaneously by a negative jump; this yields an isolated 
discontinuity at 
$x_i$ with value $\theta_U(x_i)$ that corresponds to the positive jump. 

Similarly for $\theta_L(x)$ with ``smaller value" replacing ``larger
value" in the above.

Associated with an RD touching segment is a common non-decreasing RD
section $F_T(x)\equiv F_r(x|\bar \theta(x))$ where for each
$x\in[x_1,X],\; \bar\theta(x)\equiv[\theta_L(x),\theta_U(x)].$

An {\it RD non-intersection configuration} corresponds to
the RD touchings that occur for a particular arrangement of $N$ 
disjoint touching segments
$\chi_j, (j=1,\ldots,N)$. All such segment arrangements, together with
$N=0$ for strict
(no touching) non-intersections and also $x_1=X$ for point touchings,
constitute a characterization of the entire class of RD non-intersecting 
configurations. [Not addressed is the possibility of a denumerable
number of touching segments and denumerable changes within each segment;
this might affect the completeness of the characterization.]

 An {\em analysis of intersections} appears in Appendix B.

\section{The $F(x,\theta)$ fiducial surface}\label{section4}

A three dimensional version of Figure 1 with $\theta$ presented on a
separate axis yields a surface $F(x,\theta)$ that represents the
geometric formulation. The RD $F_r(x|\theta_0)$ for given $\theta_0$
constitutes a
{\it random $\theta$ section} of the surface, i.e., the 
RD is the intersection (curve) of the surface
with the vertical plane $\theta=\theta_0$. The entire surface
$F(x,\theta)$ is generated by all the $\theta$ sections. Similarly,
each FM, $M_f(\theta|x_0)$, corresponds
to a {\it fiducial $x$ section}; all of these sections likewise
generate the surface. 
The geometric formulation in
Section \ref{section1} starts by specifying the RD $\theta$ sections,
$F_r(x|\theta)$;
these generate the surface $F(x,\theta)$ which in turn yields the 
fiducial $x$ sections. In Proposition FM of Section 1 this proves existence
of the FM $M_f(\theta|x_0)$.

To prove the geometric identity Eq.(\ref{eq3}) we note that
each individual point $(x_0,\theta_0)$ 
lies on both a $\theta_0$ section and an $x_0$ section. Hence the fiducial
measure $M_f(\theta_0|x_0)$ equals the random probability
$F_r(x_0|\theta_0)$.
Letting $\theta_0$ vary (as $\theta$) and with accompanying variable
points $(x_0,\theta)$, each fiducial $x_0$ section $M_f(\theta|x_0)$
is equal to the value at $x_0$ of
each corresponding RD $\theta$ section, $F_r(x_0|\theta),$ thereby giving
the geometric identity.

The two-dimensional {\it random plane} plot (as in Figure 1 and Figure
2(a)) of the
RD's for different $\theta$'s consists of projections of
a finite number (to be
visually meaningful) of $\theta$ sections onto the $F(=F_r)-x$
plane. They comprise superpositions of $\theta$ section RD's such as
obtained when the surface is viewed from the $\theta$ axis direction. In
the nomenclature of Subsection 3.2 this also constitutes the RD
configuration. (The geometry of the configuration is the same
whether or not the $\theta$ sections vary continuously with
$\theta$.) Similarly the {\it fiducial plane} representations, as in
Figure 2(b), are
projections of $x$-sections onto the $F(=M_f)-\theta$ plane; equivalently,
superimposed $x$-sections as viewed from the $x$ axis direction. This
constitutes an FM configuration that is dual to the RD configuration.

The plane representations provide a
clearer picture of intersecting RD's and non-monotone $m_f$'s than does the
surface representation.
The latter, however, conveys more simply the continuity requirement, with 
the three (a)-(c) continuity conditions in the Section 1 definition of
fiducial model being replaced simply by continuity of $F(x,\theta)$.

[Graves has proven in [4, Theorem 5, p.102] that for a surface
$f(x,y)$ to be continuous at a point $(x,y)$ requires that the continuity
of a $y$ section, $f(x|y)$ say, at that point be uniform in $y$. Graves
gives an
example showing that otherwise a discontinuity can result. 
It is not clear, however, whether  uniformity in $\theta$ is needed for 
the monotone distribution functions.]

\bigskip

\noindent{\bf The fiducial surface geometry}

For purpose of describing the $F(x,\theta)$ surface,
suppose that the $x$ and $\theta$ domains each consist of the unit 
interval $[0,1]$. Together with the zero-one probability range for the
ordinate $F$, the surface will then lie within the unit cube shown in
Figure 3. 


\begin{figure}[th]
\begin{center} 
\setlength{\unitlength}{0.00053333in}
\begingroup\makeatletter\ifx\SetFigFont\undefined%
\gdef\SetFigFont#1#2#3#4#5{%
  \reset@font\fontsize{#1}{#2pt}%
  \fontfamily{#3}\fontseries{#4}\fontshape{#5}%
  \selectfont}%
\fi\endgroup%
{\renewcommand{\dashlinestretch}{30}
\begin{picture}(5499,5250)(0,-10)
\path(1950,4875)(1950,1875)(4950,1875)
\path(1950,1875)(150,75)
\texture{55888888 88555555 5522a222 a2555555 55888888 88555555 552a2a2a 2a555555 
	55888888 88555555 55a222a2 22555555 55888888 88555555 552a2a2a 2a555555 
	55888888 88555555 5522a222 a2555555 55888888 88555555 552a2a2a 2a555555 
	55888888 88555555 55a222a2 22555555 55888888 88555555 552a2a2a 2a555555 }
\shade\path(1950,4275)(4350,4275)(2850,375)
	(450,375)(1950,4275)
\path(1950,4275)(4350,4275)(2850,375)
	(450,375)(1950,4275)
\path(4350,4275)(4350,1875)
\path(4350,1875)(2850,375)
\path(4350,4275)(2850,2775)
\path(1950,4275)(450,2775)
\path(450,2775)(450,375)
\path(2850,2775)(2850,375)
\path(450,2775)(2850,2775)
\dashline{60.000}(1950,1875)(3450,1875)
\dashline{60.000}(1950,1875)(450,375)
\dashline{60.000}(1950,1875)(1950,4275)
\put(300,450){\makebox(0,0)[lb]{{\SetFigFont{12}{14.4}{\rmdefault}{\mddefault}{\updefault}1}}}
\put(3000,375){\makebox(0,0)[lb]{{\SetFigFont{12}{14.4}{\rmdefault}{\mddefault}{\updefault}2}}}
\put(4425,1950){\makebox(0,0)[lb]{{\SetFigFont{12}{14.4}{\rmdefault}{\mddefault}{\updefault}3}}}
\put(4425,4275){\makebox(0,0)[lb]{{\SetFigFont{12}{14.4}{\rmdefault}{\mddefault}{\updefault}7}}}
\put(1725,4275){\makebox(0,0)[lb]{{\SetFigFont{12}{14.4}{\rmdefault}{\mddefault}{\updefault}8}}}
\put(300,2775){\makebox(0,0)[lb]{{\SetFigFont{12}{14.4}{\rmdefault}{\mddefault}{\updefault}5}}}
\put(1800,1950){\makebox(0,0)[lb]{{\SetFigFont{12}{14.4}{\rmdefault}{\mddefault}{\updefault}4}}}
\put(0,0){\makebox(0,0)[lb]{{\SetFigFont{12}{14.4}{\rmdefault}{\mddefault}{\updefault}$x$}}}
\put(5100,1875){\makebox(0,0)[lb]{{\SetFigFont{12}{14.4}{\rmdefault}{\mddefault}{\updefault}$\theta$}}}
\put(3000,2700){\makebox(0,0)[lb]{{\SetFigFont{12}{14.4}{\rmdefault}{\mddefault}{\updefault}6}}}
\put(1950,5100){\makebox(0,0)[lb]{{\SetFigFont{12}{14.4}{\rmdefault}{\mddefault}{\updefault}$F$}}}
\end{picture}
}
\caption{Inclined plane (1278) as preliminary $F(x,\theta)$ surface for FD}
\end{center} 
\end{figure}
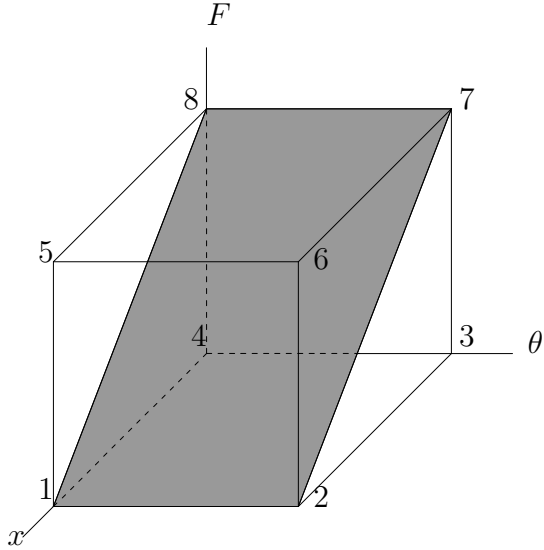

We further initially suppose the surface to be an inclined plane passing
through the vertices 1278 in Figure 3; equivalently, through the base edge
12 and its diagonally opposite edge
78. Along these two edges the RD's have probabilities of zero and one,
respectively. The RD random plane projection consists of a single 45 degree
line applicable to all $\theta$, while the fiducial plane projection
consists of horizontal parallel lines.
The RD's can be separated by pushing down side 27 toward the
base edge 23, and simultaneously pulling side 18 upward toward top
edge 58. When (approximate) coincidences with edges 23 and 58 are reached,
the RD's become complete. If the stretchings of the
inclined plane are ``smooth" (no distortions, as explained below)
the monotone random plane RD projections will be non-intersecting.
The result is a FD surface where the $x$ section $M_f (=F_f)$
curves are monotone, non-intersecting and complete.

This FD surface is symmetric about the
surface curve (with endpoints 2 and 8) that results from 
intersecting the surface with the 
vertical plane through vertical edges 26 and 48 -- equivalently, 
through the diagonally opposite vertices 2 and 8. The surface
proceeds ``monotonically'' from the base-anchored 
edges 12 and 23 to the top-anchored edges 58 and 87; vertices 2 and 8
become the symmetric
vertex anchors. The symmetric-like aspect resides in the fact that
interchanging the two surface portions  on each side of the symmetric curve
28 (which is equivalent to interchanging RD and FD) retains the monotone
character of the surface.

Opposite
monotonicity of parameter $\theta$ is obtained by changing the directions
of the $\theta$ and $x$ axes. Other variations correspond to having
symmetric base anchors at each of the four base vertices. We can state
then the following characterization of FD's: The surface for an FD
model consists of a symmetric monotone surface that progresses from two
adjacent base edges to diagonally opposite  adjacent top edges.

A non-FD surface with intersecting RD's contains depressions, or
bulges, that destroy the monotonicity. Consider, for example, the 
joined distributions in Section \ref{section2}, $F_1$ and $F_2$, with 
transition interval $\pm\theta_T$ (assuming, however, a finite domain).
Until the transition is reached one has a monotone $F_1$ surface
with section RD's $F_{r_1}(x|\theta)$ that
are monotone increasing in the interval $[\theta_m, -\theta_T]$;
the final RD is the $-\theta_T$ section $F_{r_1}(x|-\theta_T)$. After the
transition there is likewise a monotone increasing $F_2$ surface for
$\theta$ sections in the interval $[+\theta_T, \theta_M],$ starting
with $F_{r_2}(x|+\theta_T).$ In the transition interval the surface
sustains
a depression region for $x>x_T$ because of the $\theta$ change to monotone
decreasing (i.e. the $x$-partial derivative of F becomes negative) 
and continues until the terminal change at $+\theta_T$ reverts to
monotone increasing. The planar RD's intersect in the
 depression region and the planar fiducial measures are non-monotone. For 
$x<x_T$ the surface is monotone increasing.

\bigskip

\noindent{\bf Reciprocal distributions}

The fiducial surface $F(x,\theta)$ can be generalized to include
non-monotone
$\theta$ sections $m_r(x|\theta)$. In fact, a {\it reciprocal}
$m_r$ and $F_f$ relation 
can be obtained by interchanging $x$ and $\theta$ in the
operations that lead, for example, to non-FD ``composite'' random
variables, as discussed in the following Section \ref{section5}.
One obtains not only the dual FM
which is equivalent to the RD distribution, but also, after an
interchange, a valid reciprocal FM which is not equivalent to the RD's.
Confusion between these two fiducials can occur, as noted in Subsection
5.2.

\section {Composite distributions}\label{section5}

Composite variables and their distributions are the simplest and most
frequently occurring case of intersecting RD's. Moreover, their analysis
sheds
light on the nature of non-intersecting incomplete RD's, which is the
exceptional non-FD case implied by Definition CFD in Section 1. 

A simple example of a composite variable is $y\equiv |x|$: Each value of
$y$ is the  composite of $+x$ and  $-x$; the distribution of $y$ is then a
composite distribution. More generally:

\noindent {\em $y\equiv g(x)$ is a {\bf composite variable} iff $g(x)$ is
not a continuous, strictly monotone function.}

If $g(x)$ is not composite, i.e. is continuous and strictly monotone
(equivalent to 1-1, which is the generalization for several variables),
then
the distribution of $y$ is obtained from that of $x$ in the usual manner 
as $F(g^{-1}(x))$. 

We make the following conjecture:

\noindent{\em For an FD model with RD's $F_r(x|\theta)$, let $g(x)\equiv y$
and/or $h(\theta)\equiv \phi$ be composite variables. Then the composite
distributions of $y$ and/or $\phi$ comprise a non-FD model.}

We consider the following composite example, the analysis of which
indicates that the above conjecture is almost certainly true. Let 
$\theta$ be an
increasing translation parameter so that the RD and FD are given by: (For
simplicity we use the same symbol $F$ as in $F(x,\theta)$.)
\[F_r(x|\theta)=F(x+\theta) =F_f(\theta|x).\]
(More generally, with a translation parameter all subsequent definition
formulas are readily shown to apply to both RD and FM.)~~ We consider first
the composite variable
$g(x)=|x|\equiv y$, and later also $ h(\theta)= |\theta|\equiv \phi$.
We treat both the asymmetric {\it extreme value distribution (EVD)}
(writing $F_{EVD}()$ as $F()$)
 \[F(x+\theta)=1-e^{-e^{x+\theta}},\]
and also the symmetric normal distribution $N(x+\theta).$

\subsection{EVD composite I: Distribution of $y$}

The RD of $y\equiv |x|,$ denoted by $F^*_r(y|\theta)$, is obtained
from the following well known probability relations:
\begin{eqnarray}
\Pr(0\leq Y\leq y)&=& \Pr(X\in (0,y))+\Pr(X\in (-y,0))\label{eq5}\\
&=&\Pr(-y\leq X\leq y)\nonumber\\
&=& \Pr(-\infty<X\leq y) -\Pr(-\infty <X\leq -y)\nonumber.
\end{eqnarray}
Hence the composite RD and  FM are
\begin{equation}\label{eq6}
F^*_r(y|\theta)=F_x(y|\theta)-F_x(-y|\theta)=F(y+\theta)-F(-y+\theta)
=m^*_f(\theta|y).
\end{equation}
$F_x()$ denotes the distribution function when $x$ is the random variable.
(Eq.(\ref{eq6}) also yields the well known property that the (derivative)
density of $F^*_r(y|\theta)$ equals 
the sum of $F_x$ densities, $f_x(y+\theta)+f_x(-y+\theta)$.)

$F^*_r(y|\theta)$ is shown in the Figure 4(a) (undarkened) curves for
various $\theta$.
Intersections occur at every point $(y,F^*_r)$, unlike in the Section 2
joined RD's, since the transformation into
$y$ is non- 1-1 for all $x$.  The fiducial dual $m^*_f(\theta|y)$,
shown in the Figure 4(b) (undarkened) curves for various $y$, is
non-monotone.

\begin{figure} \psfig{file=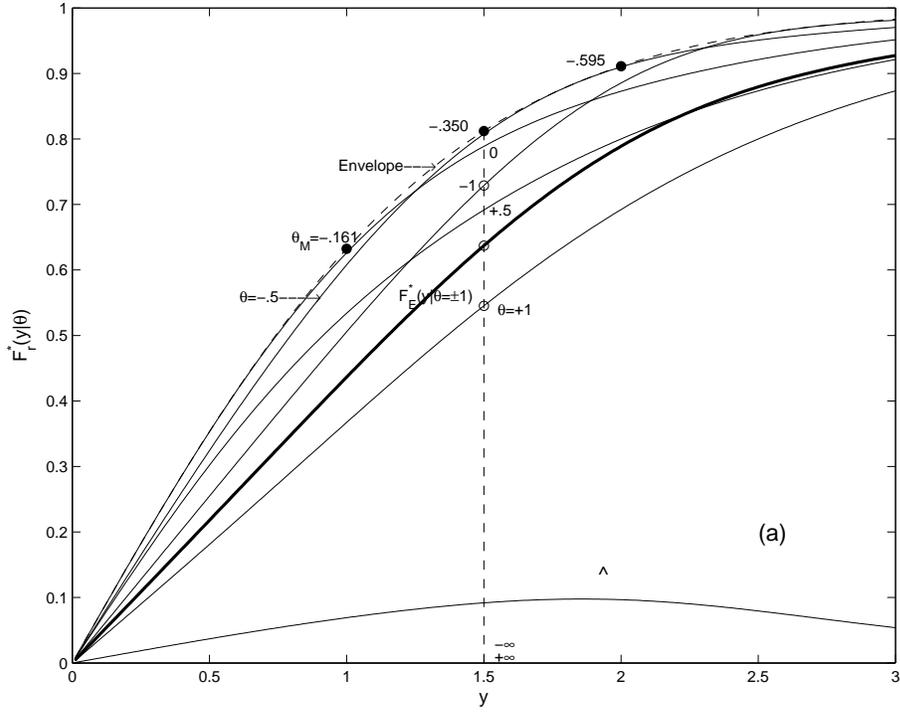,width=14.12cm,angle=0}\\
\psfig{file=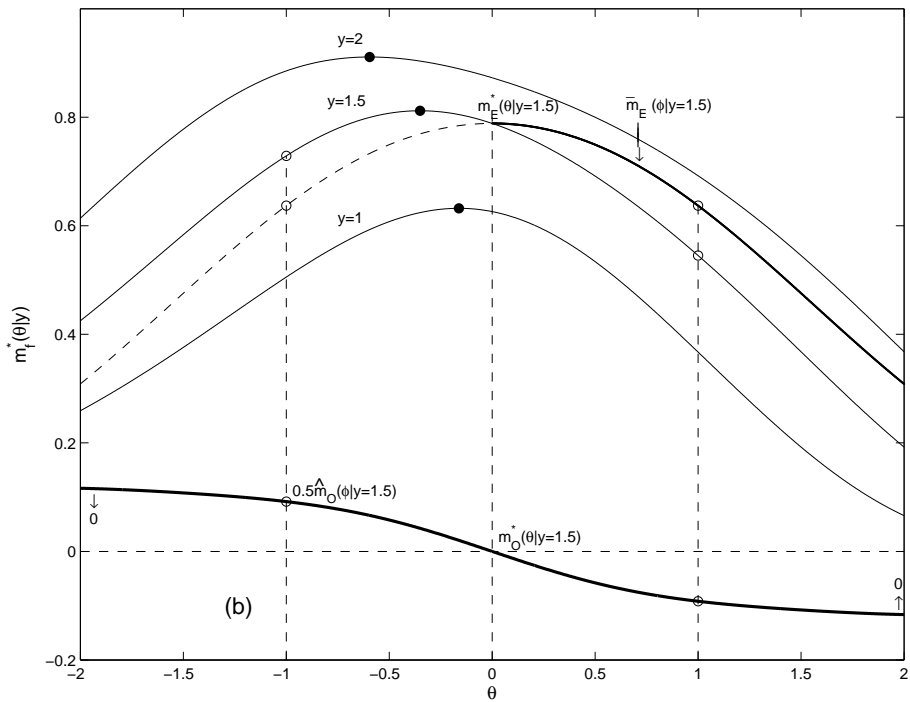,width=14.12cm,angle=0}
\caption{Composite I for EVD (a) $F_r^*(y|\theta)$ (b)
$m_f^*(\theta|y)$} 
\end{figure}

Since each intersection consists of two $\theta$ RD's, $m^*$ has for each
$y$ a maximum value, at an intermediate $\theta_M$ say, given by 
$F_r^*(y|\theta_M)(\equiv F^*_M(y))$.
In fact, equating to zero the derivative (density) of $m^*$ yields the
formula
\begin{eqnarray*}
\theta_M&=&-\ln\left({\sinh y\over y}\right)\approx
-{y^2\over6}+{y^4\over180},\\
F^*_M&=&F(y+\theta_M)-F(-y+\theta_M).
\end{eqnarray*}
The dashed curve in Figure 4(a) for $F^*_M(y)$  represents the
envelope beyond which RD probabilities will not occur. 
The values shown for
$\theta_M$ correspond to the maxima in the Figure 4(b) FM's.

$F^*_r(y|\theta)$ equals the area subtended under the EVD density 
$f(x)$ by the variable interval $(-y+\theta, y+\theta)$ with $\theta$ 
fixed, as $y$ varies from $y=0$ and $F^*_r=0$ at the EVD mode 
 ($x=0$), to $y=\infty$ and $ F^*_r=1$.  
Similarly, $m^*(\theta|y)$ is the subtended area of the same interval, 
but now with $y$ fixed and $\theta$ variable from 
$-\infty $ to $+\infty$. 
Because of the EVD left skewness, with density fall-off from the mode 
that is steeper in the positive $\theta$ direction than the negative 
$\theta$ direction, the subtended area $(-y,+y)$ at the mode is 
increased when $\theta$ is decreased until the (negative) maximum 
$\theta_M$ is reached. As $y$ is increased $\theta_M$
becomes more negative. For a symmetric density $\theta_M$ remains 
at the zero mode for all $y$.

This {\em density subtended area} procedure has various applications, for
example, in demonstrating existence of an intersection at every value of
the RD: For each $y_0>0$ and (negative) $\theta_0<\theta_M$ and with the
corresponding area of (say) $A_0 (\equiv F^*(y_0|\theta_0)$ subtended 
by $(-y_0+\theta_0, y_0+\theta_0)$, there exists a unique $\theta'_0$  with
area $A_0$ subtended by the interval 
$(-y_0+\theta'_0, y_0+\theta'_0)$.
This procedure also shows that a density with multiple modes results in
multiple $\theta$ intersections.

The dual relation between the RD and FM is illustrated by the vertical 
line in Figure 4(a) for $y=1.5$, and the corresponding $m^*(\theta)$ 
curve in Figure 4(b). Starting at $\theta=-\infty$ with $F^*_r=0$ 
(the base line in  Figure 4(a)), increasing $\theta$ yields 
increasing values of $F^*_r$ as shown for $\theta=-1$ and $-.5$,
until the maximum $F^*_M$ at $\theta_M= -.35$ is reached; 
this maximum is shown also  in the Figure 4(b)  $m^*(\theta|y=1.5)$ 
curve. Further increase in $\theta$ to 0, +1 and  $+\infty$ 
decreases $F^*_r$ and $m^*$, ultimately to zero.
\bigskip

The composite variable $\phi$ and its composite distribution will be
treated in Section 5.3 together with the related darkened curves in Figures
4(a) and 4(b), 
after analyzing the normal distribution situation.
\bigskip

\noindent {\bf Reciprocal distributions}

A {\it reciprocal distribution} is obtained by interchanging
variables $(y,\theta)$ with $(\phi,x)$. Since the FD
$F_f(\theta|x)$ is the same EVD as the RD, Eq.(\ref{eq6}) yields also a 
reciprocal fiducial composite distribution of $\phi\equiv|\theta|$, denoted
by $F^{*R}_f(\phi|x)$. This distribution is identical to the RD's in
 Figure 4(a), but now with $\phi$ as the abscissa and with $x$ replacing
$\theta$. That is,
\[F_f^{*R}(\phi|x)  =F(\phi+x)-F(-\phi+x).\]
Likewise the dual reciprocal random composite, $m^{*R}_r(x|\phi)$,
is identical to the curves in Figure 4(b), with the same interchange of
variables and hence also with the same equation above.

\subsection{Normal distribution composite}

For the normal distribution, $N(x+\theta)$, the RD $F_r^*(y|\theta)$ of the
composite $y$ 
is also given by Eq.(\ref{eq6}) with $N$ replacing $F$, and is shown in
Figure 5(a).

\begin{figure} 
\psfig{file=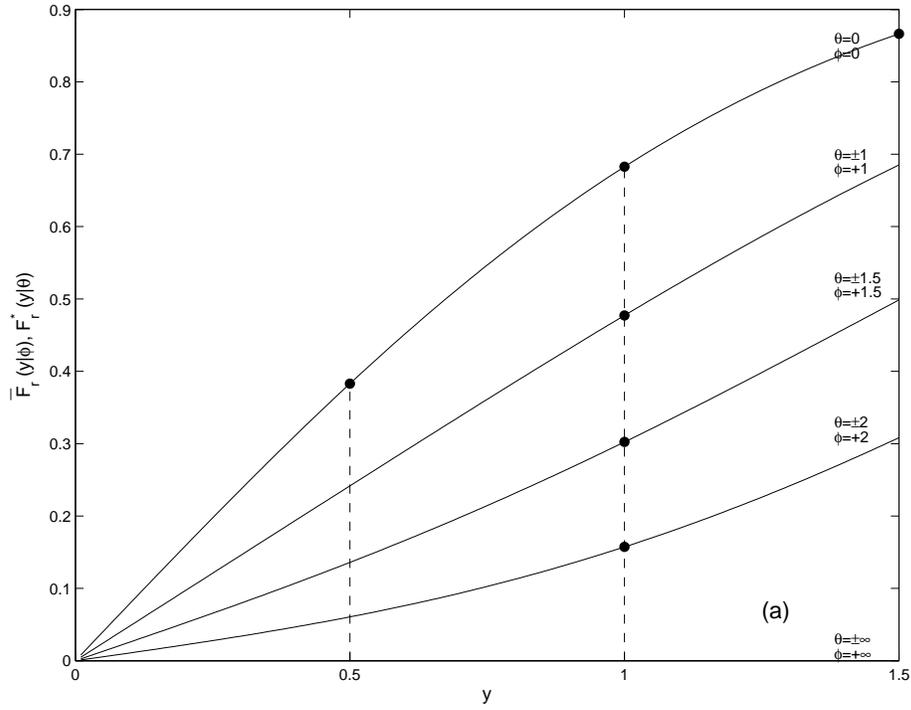,width=14.12cm,angle=0} \\
\psfig{file=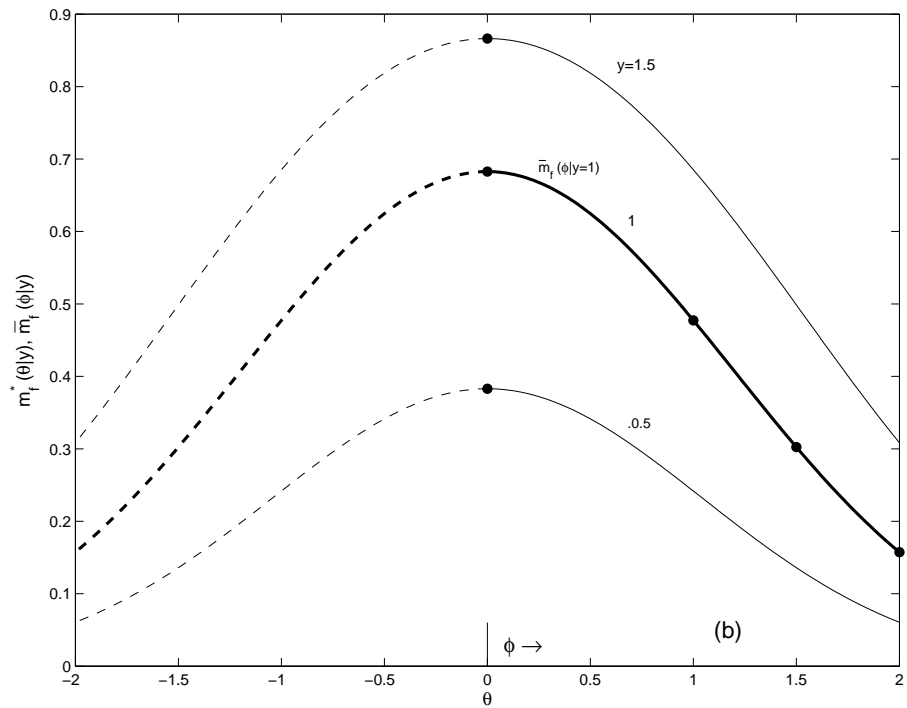,width=14.12cm,angle=0}
\caption{Composite for N.D.  (a) $F^*_r(y|\theta)$ and
$\bar F_r(y|\phi)$ (b) $m^*_f(\theta|y)$ and $\bar m_f(\phi|y)$ } 
\end{figure}

The previous EVD intersections have now come together so that the RD
for each $-\theta$  coincides with the RD for $+\theta$. (This follows from
the symmetry condition, $N(-u)=1-N(u)$, applied to Eq.(\ref{eq6}).) 
Thus these are weak intersections (defined in Appendix B), and the dual
$m^*_f(\theta|y)$ curve in Figure 5(b) 
(including both the left
dashed part and the right solid part) is non-monotone and symmetric
about the maximum at $\theta=0$; the maxima are represented also by the
(dashed) RD envelope $F_r^*(y|0)$ in Figure 5(a).

The RD coincidences in Figure 5(a) imply that the RD's can be
expressed as a function 
of the {\em composite-reducing} parameter $\phi\equiv |\theta|$.
The RD's are now  {\em composite-reduced} distributions 
$\bar F_r(y|\phi)$. They have become non-intersecting and, being
bounded by the
$\theta=0$ maximum RD $\bar F_r(\phi|0)$, they are also incomplete. 
The dual FM $\bar m(\phi|y)$ in Figure 5(b) consists of the solid
$+\theta$ portion of $m^*_f(\theta|y)$, with the abscissa parameter
$\theta$
now replaced by $\phi$; i.e. $\bar m(\phi|y)=m^*(+\phi|y),\;\phi\ge0.$
This FM is monotone, in accordance with Theorem 1, the equation being
\begin{equation}\label{eq8}
  \bar m(\phi|y)=N(\phi+y)-N(\phi-y),\ \ \ \phi,y\ge0.
  \end{equation}
The $\bar F_r$ RD's in (relabeled) Figure 5(a) decrease in magnitude
as $\phi$ increases; that is, $\phi$ is a decreasing parameter,
whether $\theta$ is initially represented as increasing or
decreasing. Hence the dual $\bar m(\phi|y)$ in (relabeled) Figure
5(b) is also monotone decreasing, from $\phi=0$ and $\bar m(0|y)(=\bar
F_r(y|0)),$ to $\phi =+\infty$ and $\bar m(+\infty|y)=0$. 
By defining $\bar \mu_f(\phi|y)\equiv 1-\bar m_f(\phi|y)$, with initial
value $\bar \mu_f(0|y) =2N(-y)$, we can obtain a
more useful representation as a monotone
increasing {\it truncated$^*$ distribution}. [The * indicates a
terminology different from that for normalized truncated distributions.
An alternative (albeit oxymoronic) terminology that emphasizes its
  relation to the RD completeness condition is ``incomplete FD''.]

The RD-FM relation can be summed up as: The composite-reduced RD's are
non-intersecting and
incomplete iff the dual FM's are monotone and
truncated$^*$.
The equivalence Eq.(\ref{eq3}) yields also: $F_r^*(y|\theta)$ is
composite-reducible
iff the dual $m_f^*(\theta|y)$ is symmetric for all $y$.

These results apply to any (density) unimodal symmetric translation
distribution
$F(x+\theta)$. For multiple modes (local density maxima), 
$\bar m(\phi|y)$ is not
monotone whence (Theorem 1) $\bar F(y|\phi)$ has intersections.

We note that $\bar m(\phi)$ is not the composite {\em distribution} of the
composite variable $\phi$, but is the FM that corresponds to the 
reduction of the RD $F_r^*(y|\theta)$, from a function of
$\theta$ to  an equivalent reduced representation as a function of $\phi$.
Appendix C shows that the composite {\em distribution} of
$\phi$ is in fact identically zero.
\bigskip

The preceding analyses, which are applicable in essence also to
non-translation
parameters, indicate the likely truth of the conjecture that
composite distributions are non-FD. Also indicated is that
reduced-composite distributions are directly related to the completeness
condition (ii) in the FD definition and in Theorem 2 of Section 3. 
Incomplete non-intersecting RD's apparently arise only from 
association with the envelope that 
exists for the non-FD RD's.
(The completeness condition (ii) in the FD definition should perhaps be
accompanied by a statement to that effect.)

Of course, a complete answer would consist of necessary and sufficient
conditions for existence of incomplete/truncated$*$ distributions.
\bigskip

\noindent
{\bf Reciprocal distribution and fiducial counterexample}

The reciprocal distribution of $\phi$ is
 obtained by starting with the FD $N(x+\theta)$ and
then applying Eq.(\ref{eq6}) as was done for $y$. The resulting
distribution is also the same and is represented in Figure 5(a) 
with variables
interchanged. The equation is:
\[F^{*R}_f(\phi|x)=N(\phi+x)-N(-\phi+x)=m^{*R}_r(x|\phi).\]

The coincidence of $+x$ and $-x$ curves in the relabeled
Figure 5(a) marks these curves as being {\it reciprocal
incomplete composite-reduced} fiducial distributions  of $\phi,\; \bar
F_f^R(\phi|y)$, with $y$ the reducing-composite variable.
The duals, $\bar m^R_r(y|\phi)$, are represented by the solid
portions of the relabeled curves in Figure 5(b).
We have then two fiducials for $\phi$:\\
(1) the dual truncated$^*$
distribution given by (with (derivative) density $ f_f()$ included):
\begin{eqnarray*}
\bar \mu_f(\phi|y)&=&1-[N(y+\phi)-N(-y+\phi)]=N(-y+\phi)+N(-y-\phi)\\
f_f(\phi|y)&=&n(-y+\phi)-n(-y-\phi)
=\sqrt{2/\pi}e^{-(y^2+\phi^2)/2}\sinh (y\phi);
\end{eqnarray*}\\
(2) the reciprocal incomplete distribution:
\begin{eqnarray*}
\bar F^R_f(\phi|y)&=&N(y+\phi)-N(y-\phi)\\
f^R_f(\phi|y)&=&n(y+\phi)+n(y-\phi)
=\sqrt{2/\pi}e^{-(y^2+\phi^2)/2}\cosh (y\phi).
\end{eqnarray*}

\noindent Since $n(y-\phi)=n(-y+\phi),$ this last equation is also the
density  $f_r(y|\phi)$ of the RD $\bar F_r(y|\phi).$

The two fiducials are relevant to Stein's counterexample [11]. 
Stein considered 
the RD of the $n$ degree-of-freedom non-central chi-square (composite)
distribution, $F_r(z|\psi)$, of $z=x_1^2+\cdots+x_n^2$ with 
non-centrality
parameter $\psi=\theta_1^2+\cdots+\theta_n^2$. For $n=1$ with  
$z=x^2, \psi=\theta^2$, the distribution is equivalent to 
the distribution of the 1-1 continuous transformation  variables 
$|x|=\sqrt z$ and $|\theta|=\sqrt \psi$, which were treated above.
Stein determined the approximate probability that the one-sided 
confidence interval for $\psi$ will cover the true $\psi$,
estimated for large $n$ from the approximating normal distribution 
with the
same mean and variance.  The probability of covering this same 
parameter was estimated also for a fiducial confidence interval 
using a similar
normal distribution approximation. The resulting large difference in
coverage
probabilities between the random and fiducial limits indicated a
serious deficiency in the latter. 

The reason for the difference was the use of the ($n$ sample version)  
reciprocal rather than the dual fiducial distribution, so that the analysis
actually confirms the non-equivalence between the reciprocal fiducial
distribution and the random distribution. Stein does remark that 
``a more reasonable fiducial distribution leads to intervals based on
... confidence intervals.''

In the absence of a fiducial framework, selecting the
reciprocal fiducial -- apparently the choice also of R.A. Fisher -- may
have been influenced by its equation being identical to the RD, a feature 
previously noted for translation-related FD's and RD's.

A criticism of Stein's probability analysis in this example has been 
given by  Pinkham [9].

\subsection{EVD Composite II: Composite distribution of $\phi$}

Previously obtained for the EVD were the asymmetrical FM's and RD's,
$m^*_f(\theta|y)(\equiv m^*(\theta)$ for brevity) and $F_r^*(y|\theta)$,
shown in Figures 4(b) and 4(a). The determination of the composite {\em
distribution} of $\phi$ from non-monotone $m^*(\theta)$ uses the fact that
an arbitrary function can be uniquely
represented as the sum of its even (E) and odd (O) component functions. For
$m^*(\theta)$ we have
\begin{eqnarray*}
m^*(\theta)&=&m^*_E(\theta) +m^*_O(\theta)\\ 
m^*_E(\theta)&=&[m^*(\theta)+m^*(-\theta)]/2\\
m^*_O(\theta)&=&[m^*(\theta)-m^*(-\theta)]/2.
\end{eqnarray*}
$m^*_E$ is symmetric; $m^*_O$ is referred to as {\em opposite}.

Plotted (darkened) in Figure 4(b) for $y=1.5$ are 
$m^*_E$ and $m^*_O$. The circular dots at $\theta=+1$ and -1 illustrate the
above equations: $m^*_E$ is the average of the $\pm\theta$ values of $m^*$,
and $m^*_O$ is one-half the difference of the two curves 
(or $=m^*(-\theta)-m^*_E(-\theta)$). The negative values for 
$m^*_O(+\theta)$ are associated with $m^*$ being a signed measure
and reflect the positive skewness of the EVD.
The above equations for $m^*_E$ and $m^*_O$ yield  
\begin{equation}\label{eq9}
m^*_E(\theta|y)=[F(y+\theta)-F(-y+\theta)+F(y-\theta)-F(-y-\theta)]/2.
\end{equation}
\begin{equation}\label{eq10}
m^*_O(\theta|y)=[F(y+\theta)-F(-y+\theta)-F(y-\theta)+F(-y-\theta)]/2.
\end{equation}
These are also the equations of the RD duals, $F^*_E(y|\theta)$
and $F^*_O(y|\theta)$.

$F^*_E(y|\theta)$, shown (darkened) for $\theta=\pm 1$ in Figure 4(a),
is the average of RD's $F^*_r(y|+\theta)$ and $F^*_r(y|-\theta).$ 
Similarly, $F^*_O$ equals one-half the difference between these two RD's
(equivalently, equals $F^*_r(y|-\theta)-F^*_E(y|\pm\theta)$)
and is shown in Figure 4(a) for $\theta=-1$; $F^*_O(y|+\theta)$ is
identical but with negative values. (Note that $F^*_O$ 
is a  non-monotone distribution.)

$F^*_E$ can also be obtained by replacing the EVD $F(x)$ by its
density $f(x)$. The replaced symmetric even component density $f_E$ is then
obtained from the even-odd
decomposition of $f$, namely, as the average of
$f(x)$ and its reflection $f(-x):\; f_E(x)=[f(x)+f(-x)]/2.$ 
The corresponding RD is 
\[F_E(x)=\int_{-\infty}^x f_E(\bar x)d\bar x=[1+F(x)-F(-x)]/2.\]
 The composite is then
\[F^*_E(y|\theta)=F_E(y+\theta)-F_E(-y+\theta)={1+F(y+\theta)-F(-y-
\theta)\over2}-{1+F(-y+\theta)-F(y-\theta)\over2},\]
which reduces to Eq.(\ref{eq9}).
Similarly, with $f_O(x)=[f(x)-f(-x)]/2$ (a signed density), 
$F^*_O(y|\theta)$ reduces to Eq.(\ref{eq10}).

As in the normal distribution, $F^*_E(y|\pm \theta)$
is a composite-reduced
distribution $\bar F_E(y|\phi)$.
The composite-reduced 
dual $\bar m_E(\phi|y)$ is shown in Figure 4(b) as the darkened positive
half of $m^*_E(\theta|y)$.

The configurations for $F^*_E(y|\pm\theta)$ and $m^*_E(\theta|y)$ for
varying $\theta$ and $y$ are, since the density $f_E$ is unimodal,
qualitatively the same as in Figures 5(a) and 5(b) for the normal
distribution.

Since $m^*$ has an $m^*_O$ component in addition to the symmetrical $m^*_E$
component,
$F_r^*$ cannot, according to a previous necessary and sufficient condition,
be represented as a composite-reduced RD.

The composite distribution, $\hat m(\phi)$, of the composite
$\phi$ can be determined from its equality to the sum of its even and odd
composites:
\begin{equation}\label{eq11}
\hat m=\hat m_E +\hat m_O.
\end{equation}

It is proven in Appendix C
that $\hat m_E(\phi)=0$ and $\hat m_O(\phi)=2m_O(+\theta)$ (or
 $=2m_O(-\theta)$ when $m_O(-\theta)$ is the positive half).
Hence  $\hat m(\phi)=\hat m_O(\phi).$

The positively valued $m^*_O(-\theta)$ in Figure 4(b) represents 
$0.5 \hat m_O(\phi)$. The dual signed RD $F^*_O(y|\theta)$, 
shown in Figure 4(a) for $\theta=-1$, represents also $.5 \hat
F_r(y|\phi)$.
\bigskip

The following summarizes the properties of the odd distribution-composite 
in  comparison with the even  reduction-composite.

\begin{itemize}
\item
$\bar m_E(\phi)=m^*(+\theta)$ and $\hat m_E=0;\;\hat m_O=2 m^*(\pm
\theta)$ and $\bar m_O$ does not exist;
\item
$F^*_E(y|+\theta)=F^*_E(y|-\theta)$ (coincident); 
$F^*_O(y|+\theta)=-F^*_O(y|-\theta)$ (opposite);
\item
$m^*_E(+\theta|y)=m^*_E(-\theta|y)$ (symmetric);
$m^*_O(+\theta|y)=-m^*_O(-\theta|y)$ (opposite);
\item
$f_E$ density symmetric; $f_O$ density opposite;
\item
For unimode $f,\; \bar F_E(y|\phi)$ is non-intersecting; 
 $\hat F_O(y|\phi)$ is
intersecting (since $\hat m^*_O(\phi)$ is non-monotone with a maximum).
\end{itemize}
\bigskip

We note finally the relation: 
\[m(+\theta)=m_E(+\theta)+m_O(+\theta)=\bar m_E(\phi)+0.5 \hat m_O(\phi).\]

\begin{center}
***************************************************************************
\end{center}

\subsubsection*{Remarks on before vs. after observation}

The previous sections have dealt with the equivalence between the 
RD and FD representations: $F_r(x|\theta)$ for all $\theta$ equals
$F_f(\theta|x)$ for all $x$. The difference in emphasis -- given $x$
vs.~given $\theta$ -- is related to their roles in applications.
When an actual observation $x_0$ occurs, the relevant information for
inference is the
single FD $F_f(\theta|x_0)$. Prior to the observation there is only the RD
probability outcomes $F_r(x|\theta)$ for all possible 
observations $x$ and for each $\theta$.

In practice this before-after  distinction may not always be sharp. In
hypothesis testing the FD representation corresponding to the
(non-observed) critical value for
acceptance is useful. The opposite may also occur, as
in confidence limits when
the appropriate $\theta$ RD can be selected after
the observation occurs.
A difference arises when the RD and FD criteria are not the same,
as in the significance test for the Behrens-Fisher problem.
It is interesting that Laplace in a 1777 paper made a 
before-after  distinction when discussing ``the search for the [proper]
mean". Laplace states (Stigler [8, p.119]): 

\begin{quote}
The problem ...  may be regarded from 
two different points of view, depending on whether we consider 
the observations before or after they are made... It is from the
first of these viewpoints that the question has, prior to now,
been treated ... However ingenious their researches have been, 
they can only be of very little use to observers.
\end{quote}

\noindent Laplace's ``observers" refers to ``geometers", such as 
astronomers, who utilize the observations. 
  
\newpage
\section*{APPENDIX A: Comparison of pivotal quantity (PQ) formulation and
geometric formulation} 

More detailed analyses 
are given in a concomitant paper on fiducial applications. A
summary appears at the end of this Appendix.

\bigskip

\subsection*{A1.  Single observation} 

It is convenient to divide PQ's into
{\em PQ1 translation parameters} and  {\em PQ2 general 
(non-translation) parameters}.

\bigskip

\noindent{\bf PQ2 general parameter formulation}

As noted by Seidenfeld [10, p.362], the PQ2 defining equation that is
considered the natural extension of PQ1 pivoting, is
\begin{equation} \label{eq41}
F(x,\theta)=\beta.
\end{equation}
In the present context $\beta$ is 
the confidence coefficient associated with a one-sided confidence limit
(either an upper limit or a lower limit). 
The PQ2  Eq.(\ref{eq41}) is
represented in Figure 1 by a horizontal line (a plane in the three
dimensional surface), say $F(\equiv F_b)=\beta$. 
This line intersects each $\theta$ RD in an $x(\equiv x_b(\theta))$ 
given by
$F_r(x_b|\theta)(\equiv\Pr(X\leq x_b|\theta))=\beta.$ 
When the RD's are strictly non-intersecting
$x_b(\theta)$ is a strictly monotone function which can be inverted to 
yield the {\em confidence limit function} $\theta_b(x)$.
In terms of random probability 
\begin{equation}\label{eq42}
\Pr(\theta_b(X)<\theta|\theta)=\beta.
\end{equation}
This equation constitutes the (non-fiducial) random probability
representation of the confidence limit function, which is analyzed in more
detail in Cram\'er [12, pp. 510-512]. Cram\'er's analysis of the inverting 
(readily modified for a one-sided limit) considers the probability
equivalence of the  sets corresponding to  $x_b(\theta)$ and $\theta_b(x)$.
These sets are related to 
the set of points $(x,\theta)_b$ of the fiducial surface that satisfy the
PQ2 Eq.(\ref{eq41}), and which are represented as  $x_b(\theta)$
in the RD $\theta$ sections and  as $\theta_b(x)$ in the fiducial $x$
sections. 
Their equality also constitutes the PQ2 formulation, represented by the
{\em PQ2 confidence limit identity}:
\begin{equation}\label{eq43}
\theta_b(x)=x_b(\theta).
\end{equation}
This fiducial-oriented equation represents the inverting in Cram\'er's
non-fiducial random-oriented treatment of confidence limit.

Letting $\beta$ vary from 0 to 1 yields the FD with
$F_f(\theta_b|x)=\beta$; that is, the $\beta$ confidence limit is
the $\beta$ percentile of the FD. Since $\theta_b(\equiv \theta(\beta))$ is
strictly monotone for each $x$, the inverse function $F_f^{-1}$ is defined
and
$\theta_b(x)=F_f^{-1}(\beta|x)$. Similarly, $F_r(x_b|\theta)=\beta$ and
$x_b(\theta)=F_r^{-1}(\beta|\theta)$. Hence the identity
Eq.(\ref{eq43}) can be written also as an {\em inverse-FD  identity} 
\begin{equation}\label{eq44}
F_f^{-1}(\beta|x)=F_r^{-1}(\beta|\theta).
\end{equation}
 This identity reflects, and is equivalent to, the fact that the PQ2
 definition
 Eq.(\ref{eq41}) relates to an ``inverse function" (say) of $F(x,\theta)$,
 namely $F^{-1}(\beta)=(x,\theta)_b$. 

Additional analysis of PQ2 is presented following the analysis of PQ1.

\bigskip

\noindent {\bf Geometric formulation} 

The FD, $F_f(\theta|x_0)$, corresponding to non-intersecting RD's is
obtained from the geometric identity Eq.(\ref{eq1}). 
The confidence limit $\theta_b(x_0)$ in Figure 1 is, as noted above, the
monotone
representation of
$(x_0,\theta)_b$ in the fiducial $x_0$ section -- essentially the same as
the inverting of $x_b(\theta)$ -- and corresponds to the $\beta$ percentile
of the FD.

For the general case of possibly intersecting RD's, we use the geometric
identity Eq.(\ref{eq3}) to define the fiducial measure $M_f$. The
confidence limit(s) are,
as before, obtained from the intersection with the horizontal plane,
$F_b=\beta$, and yield the set $(x,\theta)_b$.
For each $x$ there may be multiple $\theta$ values, say  $\theta_i(x),$
which we denote by the set $\t(\equiv\t(x))$.
Thus, we can write the geometric identity Eq.(\ref{eq3}) as
\[M_f(\t|x)=F_r(x|\t)=\beta,\]
which is interpreted as applying to each $\theta_i\in\t$.
 $\t$ then represents the ``confidence limit set''.
We consider three cases:

(1)  $\t$ is a single unique $\theta_1(\equiv\theta_1(x)).
\;M_f$ is then a strictly monotone FD $F_f$ and (iff) the RD's are strictly
non-intersecting. This is the usual situation in applications, also
assumed for the PQ2 solution.

(2)  The FD is non-strict and contains a constant interval,
$\t=[\theta_1,\theta_2]$, of admissible confidence limits.
The RD's are
touching at $x_0$ (say), and $\theta_b(x)$ is monotone but not continuous,
with a jump at $x_0$. Some (inequality) modifications are required in 
case (1) analyses.

(3)  $M_f(\equiv m_f)$ is non-monotone corresponding to intersecting RD's.
Multiple "confidence limits" are obtained, with
$\theta=\theta_1,\theta_2,\ldots$. This case
does not arise in applications. 

\bigskip

\noindent {\bf PQ1 translation parameters}

[Much of the following analysis for a strict translation parameter can be
extended to generalized translation
parameters, generated by  three strictly monotone functions.
For brevity
the general class is also referred to as translation-scale.]

The pivotal quantity  $u\equiv (x+\theta)$ was early recognized as a
significant quantity with a pivoting property that yields both the random
and the fiducial distributions.
PQ1 defines $F(x,\theta)$ in terms of $u$ and its probability distribution
$F^*(u)$:
\begin{equation}\label{eq51}
F(x,\theta)=F^*(x+\theta)\equiv F^*(u). 
\end{equation}
(When $F^*(u)$ is a probability distribution the RD's are
non-intersecting.)
This equation also defines a translation parameter as being equivalent to
existence of a PQ1 pivot.

For the RD we have
\begin{equation}\label{eq52}
F_r(x|\theta)=F^*(\x+\theta),
\end{equation}
where $\x$ 
indicates that the pivot is on $x$.
A pivoting, and with $x=x_0$, gives the FD, and also the PQ1 formulation,
\begin{equation}\label{eq53}
F_f(\theta|x_0)=F^*(x_0+\t).
\end{equation}

Note that the FD is obtained merely by interchanging
$x$ and $\theta$ in the RD. Fisher in his book [3, pp. 52-54]
illustrates the
``fiducial argument'' by using a similar interchange in an exponential
distribution with scale parameter. 

It is readily verified that the FD density of Eq. (\ref{eq53}) is identical
to the bayesian posterior
distribution that is obtained when using the non-informative uniform prior
distribution applied to the RD density of Eq. (\ref{eq52}).

The PQ1 formulation is also a geometric formulation. In
fact, the same right sides of Eq.(\ref{eq52}) (with $x=x_0$) and Eq.
(\ref{eq53}) gives the geometric identity.

The PQ1 relation for the (strict) translation parameter can be extended to
generalized translation-scale parameters.
In this family each FD is
identical to a posterior distribution with a corresponding non-informative
prior, such as $1/\theta$ for a scale parameter $\theta$.
Thus there is an (iff) equivalence relation between a PQ1 pivot, a
translation-scale parameter, and a bayesian-non-informative solution.

\bigskip

\noindent
{\bf PQ2 extension of PQ1}

It is  natural to  
extend  the translation parameter's PQ1 pivoting to
non-translation-scale parameters. This 
has been accomplished by a general PQ definition wherein the 
right side of Eq.(\ref{eq51}) is replaced by
the [0,1] uniform distribution, as in the PQ2 definition Eq.(\ref{eq41}).

The resulting properties of PQ2 differ considerably from PQ1. 
PQ1 has a real pivot with formulation identical to the geometric
formulation. PQ2's pivoting, on the other hand, is accomplished by 
monotone
inverting, which results in the FD-inverse formulation.
The PQ2 definition, which specifies the $(x,\theta)_b$ values for each
$\beta$,
is equivalent to the ``inverse'' of $F(x,\theta)$ and the resulting 
confidence limit, while
PQ1 defines a general $F^*(u)$ that ranges over all values of $F$ and
yields the FD directly. 

PQ2 is not a pivoting extension of PQ1;  equivalently,  PQ1 is not a
special case of PQ2. 
Most likely, a PQ1 extension with a pivoting property does not exist.
A similar and closely related translation-scale situation 
is the non-informative property of a
prior distribution which cannot be extended to non-translation-scale
parameters. Non-extension of pivoting is also related to the pivot in
PQ1 serving as an equivalent definition of a translation
parameter. Note, however, that the general geometric formulation is a
{\em non-pivoting} extension of PQ1's geometric formulation property.

Even more troubling is the fact that PQ2's inverse FD formulation -- which
has some merit for a single observation in being equivalent to the direct
random distribution determination of confidence limits -- is not
generalizable to 
multiple observations (treated in the next section). The situation is
analogous to defining a probability distribution in terms of its inverse,
and then attempting to  extend the inverse definition to the probability
distribution of two variables or to a conditional probability distribution.
Although perhaps not impossible, implementing any such extension would
undoubtedly be very difficult and, furthermore, unnecessary.


The principal conclusions are:  (i) The PQ2 formulation is
inappropriate and should not be used. (ii) The PQ1 pivot defines
a generalized translation-scale parameter and is also a geometric
formulation.
(iii) The geometric formulation is applicable to all
parameters and to both single and multiple observations. 

\bigskip
\subsection*{ A2. Geometric formulation derivation of FD for multiple
observations}
 
For simplicity we assume only two independent observations; the proof
for the general case is essentially the same.
Initially we suppose that the
parameters are different: $x_i$ has parameter $\theta_i$ with 
geometric identity
\begin{equation}\label{eq31}
F_{f_i}(\theta_i|x_i)=F_{r_i}(x_i|\theta_i),\\ i=1,2.
\end{equation}
The joint geometric identity and joint FD are then the product
\begin{equation}\label{eq32}
F_f(\theta_1, \theta_2|x_1,x_2)\equiv
F_{f_1}(\theta_1|x_1)\cdot F_{f_2}(\theta_2|x_2)
=F_{r_1}(x_1|\theta_1)\cdot F_{r_2}(x_2|\theta_2).
\end{equation}
That is, $\theta_1$ and $\theta_2$ are the two parameter variables 
that define the joint FD function. 

The desired FD, $F_f(\theta|x_1,x_2)$, is obtained when (i.e.,
conditional on) $\theta_1=\theta_2=\theta$ in the middle fiducial
expression in Eq.(\ref{eq32}): (For clarity the usual conditional symbol
$|$ is replaced by $||$.) 
\begin{equation}\label{eq33}
F_f(\theta|x_1,x_2)\equiv F_f(\theta_1,
\theta_2|x_1,x_2)||(\theta_1=\theta_2=\theta)=
\{F_{f_1}(\theta_1|x_1)\cdot
F_{f_2}(\theta_2|x_2)\}||(\theta_1=\theta_2=\theta).
\end{equation}
For the right RD term in Eq.(\ref{eq32}) the conditioning is
achieved
merely by substituting $\theta_i=\theta$.
Substitution does not, however, give correct conditioning for 
Eq.(\ref{eq33}). The correct result
is most readily obtained using fiducial densities.
The density equation, obtained from the derivatives of Eq.(\ref{eq32}),
merely replaces the  cumulative $F$'s in Eq.(\ref{eq33}) by density $f$'s:
\begin{equation}\label{eq34}
\{f_f(\theta|x_1,x_2)d\theta\equiv 
f_{f_1}(\theta_1|x_1)\cdot
f_{f_2}(\theta_2|x_2)d\theta_1\cdot
d\theta_2\}||(\theta_1=\theta_2=\theta).
\end{equation}

Also needed is a {\em density geometric identity} -- which is also
Fisher's identity, the density version of Eq.(\ref{eq2}) -- that 
replaces Eq.(\ref{eq31}) (with $\theta_i=\theta)$:
\begin{equation}\label{eq35}
f_{f_i}(\theta|x_i)=\frac{\partial}{\partial \theta}F_{r_i}(x_i|\theta)=
\frac{\partial}{\partial \theta} \int_{-\infty}^{x_i} f_{r_i}(\bar
x_i|\theta)d\bar x_i.
\end{equation}

The conditioning in Eq.(\ref{eq34}) is achieved, as shown below, by first
substituting $\theta_i=\theta$. The two variable $(\theta_1,\theta_2)$
density is carried into the one variable $\theta$ density. Likewise the two
variable probability density element $d\theta_1\cdot d\theta_2$ is carried
into the one variable probability density element $d\theta$. We note,
however, that the product of two densities (for the same variable) is not
in general a density, i.e. with a total integral of unity, whence the
product requires normalization by the value of the integral.  The desired
FD density should then be given by:
\begin{equation}\label{eq36}
f_f(\theta|x_1,x_2)d\theta=\frac{f_{f_1}(\theta|x_1)\cdot
f_{f_2}(\theta|x_2)d\theta}{\int_{-\infty}^{+\infty}f_{f_1}(\theta|x_1)
\cdot f_{f_2}(\theta|x_2)d\theta}.
\end{equation}

The proof of this equation uses the general conditional probability
definition applicable to arbitrary sets A, B: 
\[\Pr'(A|B)=\frac{\Pr'(A\cap B)}{\Pr'(B)}.\]
The set A consists of  arbitrarily selected pairs $(\theta_1,\theta_2)$,
 such as appears in Eq.(\ref{eq32}); A is  also the union of $A_1$
(where only $\theta_1\neq \theta_2$) and $A_2$ (where only
$\theta_1=\theta_2 =\theta)$. Set B consists of all possible values of
$\theta \in [-\infty,+\infty]$, whence $\Pr'(B)$ equals the denominator of
Eq.(\ref{eq36}). Also, since $A\cap B=A_2=\theta$ the probability (density)
$\Pr'(A\cap B)$ is the numerator of Eq.(\ref{eq36}). This completes the
proof. 

For $n$ observations the FD density consists of: (i) obtaining the $n$
fiducial densities $f_{f_i}$ from $f_{r_i}$ using the
geometric identity density Eq.(\ref{eq35}), and then (ii)
obtaining the normalized product of the $f_{f_i}$'s. 

The solution holds trivially for $n=1$ where the normalizing
denominator is unity; the solution becomes, after integrating,
the geometric identity Eq.(\ref{eq1}). 

The $\beta$ confidence limit corresponding to the solution is given by the
$\beta$ percentile of
the cumulative FD, $F_f(\theta|x_1,x_2),$ the integral of Eq.(\ref{eq36}).
The procedure is similar to that for a single observation, 
and can yield also  non-unique confidence limit solutions.
For a strict translation parameter with pivot  
 density $f^*(x+\theta)$, the fiducial solution, Eq.(\ref{eq35}) and
Eq.(\ref{eq36}), will reduce to the bayesian posterior solution 
using the non-informative uniform prior, namely, 
\begin{equation}\label{eq999}
f_f(\theta|x_1,\ldots,x_n)=\frac{\Pi f^*(x_i+\theta}{\int_\theta
\Pi f^*(x_i+\theta)d\theta}.
\end{equation}
The equivalence follows from the fact that the derivative and integral in
Eq.(\ref{eq35}) are  now inverses and hence cancel, 
yielding the defining density pivoting equation for a
strict translation parameter:
\begin{equation}\label{eq88}
f_r(x_i|\theta)=f^*(x_i+\theta)=f_f(\theta|x_i).
\end{equation}

Equivalence with the bayesian solution can be shown to hold also for
generalized translation-scale parameters.

We note the similarity between the form of the general (non-translation)
fiducial solution Eq.(\ref{eq36}), and the form of
the bayesian solution Eq.(\ref{eq999}). This suggests that
the former can be described as a {\em fiducial-bayesian} solution. In fact,
after the observations $f_{r_i}(x_i|\theta)$ have been
converted to $f_{f_i}(\theta|x_i)$ using Eq.(\ref{eq34}), the calculations
for the FD in Eq.(\ref{eq36}) are identical to the calculations for the
bayesian posterior distribution solution in Eq.(\ref{eq999}) using a
uniform prior. 

\bigskip 

\noindent{\bf Conclusions}

The main conclusions from the preceding analyses are the following:

\begin{itemize}
\item
The PQ formulation is satisfactory only for generalized translation-scale
parameters where actual pivoting occurs; the resulting fiducial
distributions are also identical to those obtained from the geometric
formulation. 

\item
For other parameters and a single observation the PQ formulation, which
entails inverting rather than pivoting, yields the {\em inverse} of 
the fiducial distribution. This property is not useful for obtaining the
fiducial distribution for multiple observations.

\item
For multiple observations, only the geometric formulation provides the
fiducial distribution and associated confidence limits, a problem
previously unsolved. For general  translation-scale parameters this
fiducial distribution
reduces to the ordinary bayesian  posterior distribution using a
non-informative prior.
\end{itemize}

The geometric formulation remedies the flaw in the PQ
formulation and  should be considered the correct fiducial formulation. 
Preliminary investigation suggests that this formulation -- besides
yielding the iff result on existence of fiducial distributions and  also
the solution for multiple observations -- can be expected to yield
additional new results. This may well lead to increased acceptance of
fiducial probability
and recognition of its essential role in statistical theory and practice.

\bigskip

\begin{center}
***************************************************************************
\end{center}

\vbox{
\noindent{\bf Fiducial applications paper}

A related concomitant paper presents an overview of selected fiducial 
applications using the geometric formulation.
Besides providing a more  detailed treatment of confidence limits and FD's,
other areas, especially significance tests, are addressed. In particular, 
it is shown that a null hypothesis is accepted if and only if the 
hypothesis lies within the upper (and/or lower) confidence limit.
For sequential selection between two alternative hypotheses,
an additional ``no-difference'' hypothesis is required; termination occurs
when one of the three hypotheses satisfies the upper and/or
lower confidence limit (the reverse of the limit for null hypotheses).
}

Other topics include a detailed treatment of PQ's, generalized translation
parameters and bayesian
solutions, and analysis of adjacent FD's. 

Application  of the geometric formulation to multiple parameters is
expected to address the following topics: extension of the non-intersecting
FD existence
theorem; FD's and confidence regions  for single and joint  parameters;
extension  of Lindley's characterization of bayesian parameters; and
optimality of significance tests. Also to be considered is the problem of
non-unique PQ solutions that have been addressed especially by Savage and
Tukey in [13]. 

\section*{APPENDIX B: Analysis of intersections}

The proof of Theorem 1 in Section 3 yields also the following
contrapositive theorem:
\bigskip

\noindent{\bf Theorem 1CP} {\em
$M_f(\theta|x)$ is non-monotone for some $x(=x_0)$ say, iff
$F_r(x|\theta)$ has an intersection for some $x'_0\;( =$ or $\neq x_0).$}

A better understanding of intersections can be gained from an
independent proof of this theorem. Intersections, however, are more
varied and complex
than non-intersections and a complete analysis would exceed the scope of
this paper. 
The preliminary analysis here will, however, lead to a proof of the 
above theorem.

It is useful to note the following three cases pertaining to right and left
neighborhoods of $x_0$, for two $\theta$ RD's which
have $x_0$ as a common point, i.e., which satisfy condition (i):
$F_r(x_0|\theta_1)=F_r(x_0|\theta_2)$. Condition (ii) for each case is:

\noindent
(ii') touching: the same $\theta$ inequality in left and right
neighborhoods; touching then applies to all $\theta$ in the interval
$[\theta_1,\theta_2].$\\
(ii'') ordinary intersection: opposite $\theta$ inequalities in left
and right neighborhoods.\\
({{{ii}'}'}') weak intersection: equality in both left and right
neighborhoods: $F_r({x-}|\theta_1)= F_r({x-}|\theta_2)$ 
and similarly for ${x+}$.\\
We also consider the mixed combination with equality 
({{{ii}'}'}') on one
neighborhood side and inequality (ii'') on the other side.

The non-monotone (NM) contrapositive of the general monotone definition is:
\smallskip

\noindent{\bf Definition NM}  {\em $f(z)$ is a non-monotone 
function if there exist $z_i\; (i=1,2,3,4)$ (not necessarily all different)
such
that  $z_1>z_2$ implies (say) $f(z_1)>f(z_2)$ and $z_3>z_4$ implies
$f(z_3)<f(z_4)$.}

The following theorem may serve as a more useful equivalent
definition:

\noindent{\bf Theorem EQ/NM} {\em A continuous function $f(z)$ is
non-monotone iff there
exist $z'$ and $z''$ such that $f(z')=f(z'')$, provided that
$f(z)$ is not constant on the interval $[z',z''].$}

The proof of the theorem when $f(z)$ is non-monotone is essentially 
the same as in Theorem 1S where Lemma 1 is cited together with
the various triplet cases of the $f(z_i)$. In the other direction 
(without continuity), there exists $z_0$ such that (say) $z'<z_0>z''$; then
the Definition NM conditions are satisfied for the pairs $(z',z_0)$ and
$(z_0,z'')$.


\noindent{\bf Ordinary intersection}

We have the following

\noindent{\bf Proposition 1}\\
{\em $M_f(\theta|x_0)$ is non-monotone if a finite number 
$n(\geq2)$ of RD's intersect at $x_0;$ one may then take $x'_0=x_0$ in
Theorem 1CP.}

\noindent
In fact, we have $F_0=F_r(x_0|\theta_i)(\equiv M_f(\theta_i|x_0)),$
$i=1,\ldots,n$ but not for any $\theta^*\ne\mbox{ all }\theta_i.$ Theorem
EQ/NM then applies.

In the Section \ref{section2} example the point labeled $A$ in
Figure 2(a) has three RD
intersections at $x_0=1.25$; the corresponding $\theta_1, \theta_2,
\theta_3$ are shown in the non-monotone FM of Figure 2(b).
It is unnecessary to verify the intersection reverse
inequalities (ii'') in the right and left neighborhoods; if the
inequalities were not reversed the RD's would be touching. 
\bigskip

\noindent{\bf Proper interval intersection}

Suppose that the number of intersections at $x_0$ 
is infinite, consisting of all $\theta$'s in a $\theta$ 
interval $[\theta_1,\theta_2]$, which is assumed to be a proper subset of
$[\theta_m,\theta_M]$. Since $M_f$ is constant in the 
$\theta$ interval, we have the exceptional case in Theorem EQ/NM:
$M_f$ can be monotone
when the two remaining segments, $[\theta_m,\theta_1]$ and
$[\theta_2,\theta_M]$, are monotone (increasing, say). 

This situation arises in the
Section \ref{section2} example when $x_0 =x_T$, the boundary
transition point between intersections and non-intersections (see
Figure 2(a)). All transition RD's in
the intersection region intersect at $x_T$, and the remaining
RD's are increasing. The intersection inequality, together with
the continuous variation with $x$ of the RD's, implies that the left
(say) neighborhood of $x_T$ will be intersection-free with monotone
$M_f$'s, while the right neighborhood will have intersections at any
$x'_0 (> x_0)$ from RD's within the intersection
region. These satisfy Proposition 1 and thus $M_f$ is non-monotone.

Essentially the same argument applies when either, but not both,
$\theta_1=\theta_m$ or $\theta_2=\theta_M.$

The above remarks yield
\smallskip
  
\noindent{\bf Proposition 2}\\
{\em (a) If the intersections at $x_0$ do not consist of a single interval
then $M_f(\theta|x_0)$ is non-monotone.\\
(b) If the intersections at $x_0$ consist of a single interval
then there is an $x'_0$ in the right or left neighborhood of $x_0$ 
where $M_f(\theta|x'_0)$ is non-monotone.} 

\noindent  Equivalently, (a) is `iff' whence (b), also `iff', is the
contrapositive of (a). 
\bigskip

\noindent{\bf Complete interval intersection}

When all $\theta$ RD's in the complete
(finite) interval $[\theta_m,\theta_M]$ intersect at $x_0$, then 
$F_r(x_0|\theta)=$constant$=M_f(\theta|x_0)$ for all $\theta$
and hence are trivially monotone (with unspecified $\theta$ direction).
Since no RD's remain to provide intersections in a neighborhood,
monotones in $\theta$ can exist, in opposite directions, 
for each side of
$x_0$. Thus the intersection at $x_0$ may not lead to a non-monotone 
$M_f$ for any $x'_0\neq x_0$, seemingly counter to Theorem 1 
(or Theorem 1CP). 

As noted previously, this {\it endpoint solution} may in fact arise, where 
$x_0$ becomes a common endpoint, after
normalization, for two adjacent FD's. 
However, the solution is not applicable to Theorem 1 (or Theorem 1CP) 
which implicitly requires a single fiducial
model with zero-one RD endpoints.
\bigskip

\noindent{\bf Weak intersection} 

Suppose the $\theta_1$ and $\theta_2$ RD's satisfy an
equality $F_r(x|\theta_1) =F_r(x|\theta_2)$ not only for  
${x-}$ and ${x+}$ in left and right neighborhoods of $x_0$,
but also when these neighborhoods encompass the entire domain. The
$\theta_1$ and $\theta_2$ RD's then completely coincide for all
$x$. Since
$M_f(\theta_1|x)= M_f(\theta_2|x)$, it follows that $M_f$ is 
non-monotone for all $x$ and the RD's have a so-called {\it  weak
intersection}. 

An example appears in Section 5, Figures 5(a) and 5(b), where weak
intersection coincident RD pairs are obtained for each $\pm \theta$.

Weak intersection coincidence can be extended to a finite number of
$\theta$'s. However, an extension to a $\theta$ interval, with
accompanying constant FM, yields a non-intersection touching interval. 

In view of Proposition 2,
Proposition 1 can now be strengthened to

\bigskip

 \noindent {\bf Proposition $1^*$}   {\em
$M_f(\theta|x_0)$ is non-monotone iff there are ordinary or weak RD
intersections at $x_0$.} 

 \noindent (This proposition assumes that it is not admissible to
have two disjoint
interval intersections at $x_0$, or a single interval intersection plus 
discrete intersections.)

An {\it intermediate weak intersection} has a 
restricted coincidence, i.e. the neighborhood
equalities hold only for the interior of a limited intersection 
$x$-segment, $[x_1,x_2]$; also the inequality
in the left neighborhood of $x_1$ is the opposite of the inequality in
the right neighborhood of $x_2.$ At the endpoints will be an ordinary
intersection while at interior $x$'s one has weakly intersecting RD's.

\bigskip

 A number of areas require further investigation, including $\theta$ 
additions and
removals in an $x$-segment; also whether 
intersecting $\theta$ RD's can be replaced by touching $\theta$-intervals.
The various uncertainties  associated with
intersections do not, however, affect the following
proof of Theorem 1CP, which essentially combines the previous
results. 

\bigskip

\begin{proof1CP}
Suppose there are intersections at  $x_0$ with
$F_r(x_0|\theta)=F_0 (=M_f(\theta|x_0$)) for all  $\theta$ 
in a closed set $\Theta$. (a) If $\Theta$ is not an interval then
there
exist at least two $\theta'$s, say $ \theta_1, \theta_2$, such that
$M_f(\theta_i|x_0)=F_0\;(i=1,2)$, and also other $\theta$'s in the interval
$[\theta_1,\theta_2]$ that are not equal to $F_0.$  Theorem  EQ/NM then
states that $M_f$ is non-monotone. (b) If $\Theta$
is an interval (and a proper subset of $[\theta_m,\theta_M])$ with
$M_f(\theta|x_0)=$constant for all $\theta$ in $\Theta$, one has a
(proper) $\theta$-interval intersection. In a 
neighborhood of $x_0$ there will then exist an $x'_0$ with 
non-interval intersections and non-monotone $M_f(\theta|x'_0)$.

Conversely, suppose there exists an $x_0$ and an  FM
such that $M_f(\theta|x_0)$ is non-monotone. Then there exist 
at least two (non-touching) $\theta_i$ such that $M_f(\theta_1|x_0)=
M_f(\theta_2|x_0)$ and
not contained in an intersecting $\theta$-interval. Hence the $\theta_i$
RD's must be intersecting. 
\end{proof1CP}

\section*{APPENDIX C:  Proof for composite distributions}

This appendix derives the distributions for  $\hat m_E$ and $\hat m_O$ in
Eq.(\ref{eq11}). 
The proof 
uses the theorem, proven in [5, Theorem 4.4, p. 231], that a non-monotone
function, such as $m^*_E(\theta')$ in Figure 4(b),
can be represented as the difference of
two monotone functions, ${\cal M}_1(\theta')$ and ${\cal M}_2(\theta')$
with ${\cal M}_1(\theta') - {\cal M}_2(\theta')=m^*_E(\theta')$.
(In the following it is convenient to let 
$\theta$ have positive values only, with negative values denoted
by $-\theta)$; otherwise (as in the above equation) $\theta'$ is used. In
some equations either $\theta$ or $\theta'$ is acceptable.)

An example of the construction, such as for $m^*_E(\theta')$, 
adjoins to the increasing negative half $m^*_E(-\theta)$ at its 
maximum value at $\theta=0$ --  equal, say, to $A\equiv m^*_{E}(0)$ -- 
the reflection of the decreasing half $m^*_E(+\theta)$. This defines the
monotone 
function ${\cal M}_1(\theta')$ =  $2A-m^*_E(+\theta)$ for  $\theta'>0$,
and $=m^*_E(-\theta)$ for $\theta'<0$. 
We define also the monotone function 
${\cal M}_2(\theta')=2A-2m^*_E(+\theta)$  for $\theta'>0$ and = 0 for 
$\theta'<0$. Then 
${\cal M}_1(\theta')-{\cal M}_2(\theta') =m^*_E(\theta')$, as required. 

Applying Eq.(\ref{eq5}) to ${\cal M}_1$ and ${\cal M}_2$, the difference
then yields $\hat m_E=\hat{\cal M}_1-\hat{\cal M}_2$. In fact, the
right and left probabilities about zero are
$R_i={\cal M}_i(+\theta)-{\cal M}_i(0)$ and
$L_i={\cal M}_i(0)-{\cal M}_i(-\theta)$, which gives then
$\hat{\cal M}_i=R_i+L_i$. 

For $m^*_E$ we have   
\begin{eqnarray*}
&&R_1=A-m^*_E(+\theta),\; L_1=A-m^*_E(-\theta),\; \hat {\cal
M}_1=2A-2m^*_E(\theta),\\
&&R_2=2A-2m^*_E(+\theta),\;L_2=0,\;\hat {\cal M}_2=2A-2m^*_E(\theta),
\end{eqnarray*}
whence $\hat m^*_E=[{\cal M}_1-{\cal M}_2](\theta)=0.$

The general bounded variation construction in [5] for the ${\cal M}_i$
is specialized here for a continuous (non-monotone) function 
$m(\theta')$  with alternating increasing and
decreasing segments. The monotone increasing construction for
${\cal M}_1$ is obtained by reflecting
each monotone decreasing segment, which is then joined 
monotonically to the preceding increasing segment.
${\cal M}_2$ is constant at the increasing segments and  equal in 
magnitude to twice each reflected
decreasing segment (as in the above $m^*_E$ example), and with all
segments  joined monotonically. The result is ${\cal M}_1-{\cal M}_2
=m(\theta')$.

This construction has a signed measure interpretation where
the positive
measure ${\cal P}$ corresponds to increasing $m$ and negative
measure ${\cal N}$ corresponds to decreasing $m$, with the change 
occurring at each local maximum or minimum. Then 
${\cal M}_1={\cal P}+{\cal N}$ and ${\cal M}_2=2{\cal N}$, 
yielding the signed measure
representation $m={\cal M}_1-{\cal M}_2 ={\cal P}-{\cal N}$. 

Additivity of the 
signed measure and the ${\cal M}_i$  enables one to alter the sequence of
operations. Instead of $R_i+L_i=\hat {\cal M}_i$ and $\hat m=\hat
{\cal M}_1-\hat {\cal M}_2$, as in the above equations, one can use:
\begin{eqnarray*}
&&R_{1,2}\equiv R_1-R_2=[{\cal M}_1-A]-[{\cal M}_2-0]=[{\cal M}_1-{\cal
M}_2](+\theta)-A=m_E^*(+\theta)-A,\\
&&L_{1,2}\equiv L_1-L_2=[A-{\cal M}_1]-[0-{\cal M}_2]=
A-[{\cal M}_1-{\cal M}_2](-\theta)=A-m_E^*(-\theta),
\end{eqnarray*}
\[\hat m_E=R_{1,2}+L_{1,2}=0.\]

Thus one need not obtain the actual equations for the ${\cal M}_i$,
since the difference always equals the value of (non-monotone)
$m(\theta')$. Note also the canceling of the non-zero value A.
We remark also that the analysis is unaffected by  an ordinate translation
or  reflection.

For the analysis of the odd function $m^*_O$ in Figure 4(b) it is
somewhat simpler to treat the reflected function $m_O(\theta')\equiv
m^*_O(-\theta')$, so that the  positive values appear for $+\theta$.
To  obtain ${\cal M}_1$, the two decreasing segments are reflected: the 
segment $\theta'>\theta_M$ (where the $m_O(+\theta)$ maximum occurs) and
the $\theta'<-\theta_M$ segment. For the unreflected monotone increasing
segment, 
$-\theta_M<\theta'<+\theta_M,$ we have ${\cal M}_1=m_O$ and also ${\cal
M}_2(\theta') =0$. ${\cal M}_2$ is otherwise increasing, in a
manner that ensures $[{\cal M}_1-{\cal M}_2](\theta')=m_O(\theta')$. Then
\[R_{1,2}(+\theta)=[{\cal M}_1-{\cal M}_2](+\theta)=m_O(+\theta).\]
Since $L_i[{\cal M}_i(-\theta)]=0-{\cal M}_i(-\theta)$ we have also
\[L_{1,2}(-\theta)=-[{\cal M}_1-{\cal M}_2](-\theta)=-m_O(-\theta)
=m_O(+\theta),\]
\[\hat m_O(\theta)=2 m_O(+\theta).\]
For the unreflected $m^*_O$, the right and left measures are reversed and
$\hat m^*_O(\theta')=2 m^*_O(-\theta)$.
This completes the composite solutions for Eq.(\ref{eq11}).

$\hat m$ can also be obtained from a direct construction of ${\cal M}_1$
and ${\cal M}_2$. When the maximum occurs for negative
$\theta$ -- as in $m^*$ of Figure 4(a) -- the composite is given by 
$\hat m(\theta)=m(-\theta) -m(+\theta)$; by definition of odd function
this equals  $2m_O(-\theta)=\hat m_O.$
\bigskip

A complete solution for an arbitrary even function requires also
consideration of a function, $m_{E^o}$ say, with $m_{E^o}(0)=
m_{E^o}(\pm\infty)$. The  simplest such  function can be obtained
from the previous odd function
by reflecting $m_O(-\theta)$, the negative half of $m_O$, while
retaining the positive half. That is,
\[m_{E^o}(+\theta)=m_{E^o}(-\theta)=m_O(+\theta).\]
The construction for $m_{E^o}$ differs from the previous monotone 
$m^*_E(+\theta).\; {\cal M}_1[m_{E^o}(\theta')]$ is the same as in
$m_O$, as is also ${\cal M}_2(+\theta)$, whence $R_{1,2}$ remains the
same, equal to $m_{E^o}(+\theta)$. For left measures, even though the
increasing and decreasing segments for $m_{E^o}(-\theta)$ 
are the reverse
of $m_O(-\theta)$, we still have
\[L_{1,2}(-\theta)=-[{\cal M}_1-{\cal M}_2](-\theta)=-m_{E^o}(-\theta),\]
 Hence, as for
$m^*_E$, we have $\hat m_{E^o}=0$.
Thus the composite distribution for any symmetric function is zero.

The question arises whether 
a symmetric function $m_E(\theta|y)$ with $m_E(0)=0$, similar to $m_{E^o}$,
can be realized (with varying $y$) as a  translation-related
composite FM. Using the ``density subtended area'' procedure with the 
varying $\theta$ interval $[-\theta+y,+\theta+y]$, one can readily verify
that the symmetric RD density defined by $f(x)=0$ for 
$-a<x<+a$ and otherwise non-zero symmetric, yields the following FM's: For
$y<a,\; m_E(\theta'|y)=0$ for $|\theta'|<a$ (a modified $m_{E^o}$); 
for $y=a,\; m_E=m_{E^o}$; and for
$y>a,\;m_E$ is similar to $m^*_E$ in Figure 4(b).

The previous analysis for $m^*_O$ applies also to odd functions where
the limiting values for plus infinity and minus infinity are not 
the same. (This case represents a 
non-monotone generalization of the ordinary application of 
Eq.(\ref{eq5}) in obtaining the distribution of $y$.)
The composite distribution for any odd function $m_O(\theta')$ is then also
given by $2m_O(+\theta).$
\bigskip

\section*{References}

\noindent [1] Cram\'er, H. (1946), {\it Mathematical Methods of
Statistics}, Princeton University Press.

\noindent [2] Fisher, R.A. (1930),  ``Inverse probability," 
{\it Proc. of the Cambridge Phil. Society,} 26,  528-535.

\noindent [3] Fisher, R.A. (1956), {\it Statistical Methods and Scientific
Inference}, Edinburgh: Oliver and Boyd.

\noindent [4] Graves, L.M. (1946), {\it The Theory of Functions of Real
Variables}, First Edition, New York, McGraw Hill.

\noindent [5] Haaser, N.B. and Sullivan, J.A. (1971), {\it Real Analysis}, 
Van Nostrand.

\noindent [6] Halmos, P.R. (1950),  {\it Measure Theory}, New York, Van
Nostrand.

\noindent [7] Hannig, J. (2009), ``On generalized fiducial
inference," {\it Statistica Sinica,} 19, 491-544.

\noindent [8] Lindley, D.V. (1958), ``Fiducial distributions and
Bayes' theorem," {\it J. of Royal Stat. Soc., Series B} 20, 102-107.

\noindent [9] Pinkham, R.S. (1966), ``On a Fiducial Example of C. Stein,
{\it J. of Royal Stat. Soc., Series B} 28, 53-54.

\noindent [10] Seidenfeld, T. (1992), ``R.A. Fisher's fiducial argument and
Bayes' theorem," {\it Stat. Sci.} 7, 358-368.

\noindent [11] Stein, C. (1959), ``An example of wide discrepancy between
fiducial and confidence intervals," {\it Ann. Math. Stat.,} 30, 877-880.

\noindent [12] Stigler, S.M. (1990), {\it The History of Statistics: The
Measurement of Uncertainty before 1900}, Harvard University Press.

\noindent [13] Tukey, J.W. (1957), ``Some examples with fiducial
relevance," {\it Ann. Math. Stat.,} 28, 687-695.

\subsection*{Acknowledgment}
My thanks and appreciation go to Professor Refael Hassin. Rafi
generously took upon himself various responsibilities, including
computerizing the figures and converting many drafts, which used
antiquated
equation software, into error-free pdf files. Without his encouragement
and patience, his numerous technical and non-technical suggestions, and
his excellent advice to not wait until the concomitant applications paper
was completed, this present paper would probably not exist, now nor in the
foreseeable future.

\end{document}